%% file: main.tex
\documentclass[11pt]{amsart}
\usepackage[utf8]{inputenc}
\usepackage{amssymb,amsmath}
\usepackage{graphicx}
\usepackage{color}
\usepackage{tikz}
\usetikzlibrary{matrix}
\usepackage{pgfplots}
\usepackage{enumitem}
\usepackage{setspace}
\usepackage{algpseudocode}
\usepackage{mathtools}
\usepackage[T1]{fontenc}
\usepackage{xfrac}

\usepackage{accents}
\usepackage{multicol} 
\usepackage{multirow}
\usepackage{soul}
\usepackage{subcaption}
\usepackage{booktabs}
\usepackage{array}
\usepackage{comment}
\usepackage{algorithm2e} 
\newcolumntype{L}[1]{>{\raggedright\let\newline\\\arraybackslash\hspace{0pt}}m{#1}}
\newcolumntype{C}[1]{>{\centering\let\newline\\\arraybackslash\hspace{0pt}}m{#1}}
\newcolumntype{R}[1]{>{\raggedleft\let\newline\\\arraybackslash\hspace{0pt}}m{#1}}
\usepackage{siunitx}
\usepackage{bbold}

\usepackage[hidelinks]{hyperref}
\usepackage[nameinlink,noabbrev]{cleveref}
\newtheorem{theorem}{Theorem}

\theoremstyle{definition}

\theoremstyle{lemma}

\theoremstyle{remark}
\newtheorem{remark}[theorem]{Remark}

\Crefname{assumption}{Assumption}{Assumptions}
\numberwithin{theorem}{section}
\numberwithin{equation}{section}
\numberwithin{table}{section}
\numberwithin{figure}{section}
\input{def}
\AddToHook{cmd/appendix/before}{\crefalias{section}{appendix}}
\definecolor{corrRed}{RGB}{18,124,175} 

\allowdisplaybreaks
\begin{document}
\title[Three-term iterations for energy-based models]{Three-term Recurrence Iterations\\ for Energy-based Models}
\author[]{R.~Altmann$^\dagger$, J.~Ramme$^\ddagger$, P.~Schulze$^\ddagger$}
\address{${}^{\dagger}$ Institute of Analysis and Numerics, Otto von Guericke University Magdeburg, Universit\"atsplatz 2, 39106 Magdeburg, Germany}
\address{${}^{\ddagger}$ Institute of Mathematics, Technische Universität Berlin, Str.~des 17. Juni~136, 10623 Berlin, Germany}
\email{robert.altmann@ovgu.de, ramme@math.tu-berlin.de, pschulze@math.tu-berlin.de}
%
\date{\today}
%
%
\begin{abstract}
It is well-known that the midpoint rule preserves the dissipation inequality if applied to a certain class of energy-based models. We introduce an appropriate scaling of the state variables such that the symmetric part of the resulting iteration matrix is guaranteed to be positive definite. This allows the application of three-term iteration schemes such as the methods of Widlund and Rapoport. Special emphasis is put on examples where the symmetric part is block diagonal such that the computations decouple. This then leads to efficient dissipation-preserving numerical schemes as illustrated in two numerical examples, namely the biharmonic heat equation and linear poroelasticity. 
\end{abstract}
%
%
\maketitle
%
%
%
%
%
\section{Introduction}
	
The numerical simulation of dynamical systems is employed in various application areas to predict the behavior of physical or other systems.
Such systems are typically modeled as ordinary differential equations or, more general, differential--algebraic equations (DAEs) where the algebraic constraints may arise, e.g., from constitutive relations, coupling equations, or boundary conditions~\cite{MehZ24}.
After applying an implicit time discretization scheme to such a system, solving the time-discrete system typically requires the solution of a linear equation in each time step or, in the case of nonlinear problems, in each iteration of a Newton-type solver.
Especially when the dynamical system arises from the spatial discretization of a partial differential(--algebraic) equation, the dimension of these linear equation systems may be very large.
In such cases, it is often necessary to use iterative solvers, since they save both computational effort and memory usage.
Another motivation for iterative schemes is that their termination criteria are based on user-specified tolerances, which may be exploited to balance the algebraic error of solving the linear system with the discretization error~\cite{AriLMS13}.
A popular class of iterative solvers are Krylov subspace methods; see~\cite{Saa03,LieS12} for a detailed overview. In this context, it is desirable to obtain a short recurrence formula for the update of the solution to significantly limit the memory usage and computational effort.

For linear systems $\Amat x = b$ with a symmetric matrix $\Amat$, 
two well-known examples are the minimal residual method (Minres) and, when $\Amat$ is also positive definite, the conjugate gradient (CG) method. Both methods rely on three-term recurrences while simultaneously minimizing the residual and the error over a (shifted) Krylov subspace, respectively.
For non-symmetric matrices $\Amat=\Hmat+\Smat$ with positive definite symmetric part $\Hmat=\tfrac12(\Amat+\Amat^T)$ and skew-symmetric part $\Smat = \tfrac12(\Amat-\Amat^T)$ similar three-term recurrence methods have been proposed in \cite{ConG76,Rap78,Wid78}. 
The method introduced in \cite{ConG76,Wid78}, which we will refer to as the method of Widlund, yields a three-term recurrence and an error minimization property of the iterates as demonstrated in~\cite{Eis83}.
A variational analysis of this method and other CG-like techniques is provided in~\cite{SzyW93}.
An alternative approach was proposed by Rapoport in \cite{Rap78}. In contrast to the method of Widlund, this approach is based on a residual minimization property.
If $\Hmat$ is not necessarily positive definite but $\Smat$ has small rank, another Krylov subspace method has been proposed in \cite{BecR08}. This method is based on residual minimization and yields a short recurrence as well. 

In~\cite{ManM21}, the authors consider a system where $\Hmat$ is only assumed to be positive semi-definite. 
They introduce a two-level iterative scheme based on preconditioning from the left and from the right, which results in a shifted skew-symmetric matrix for which methods with short recurrences and residual minimization are available.
Another approach for the case with positive semi-definite $\Hmat$ is presented in \cite{GueLMS22} where the authors first transform the matrix pencil $(\Hmat,\Smat)$ to a staircase form and then perform a block Schur complement reduction.
This results in a decoupled system where the nonsingular part of $\Hmat$ can be treated using one of the methods for positive definite $\Hmat$ and the other blocks may be treated by methods for purely symmetric positive definite or purely skew-symmetric problems. 
While the methods of Widlund and Rapoport are based on preconditioning from the left, the authors in~\cite{DiaFK23} present an approach based on preconditioning from the right for systems with positive definite symmetric part.
As for the methods of Widlund and Rapoport, every iteration in the method of~\cite{DiaFK23} requires linear system solves with $\Hmat$ for which the authors propose inexact solves, e.g., via an incomplete Cholesky decomposition of $\Hmat$ or an inner CG-iteration. This then yields flexible variants of the methods of Widlund and Rapoport.
Recent works consider the extension of the methods of Widlund and Rapoport to the infinite-dimensional case \cite{MehSS25ppt} and a combination of Widlund's method with deflation techniques to increase the convergence speed \cite{DuFW25}.

In~\cite{ManM21,GueLMS22,DiaFK23,MehSS25ppt}, the authors motivate their methods by the fact that systems with positive (semi-)definite $\Hmat$ arise in the time discretization of dissipative or port-Hamiltonian systems.
These are dynamical systems with a special algebraic or geometric structure which guarantees an energy balance, see, e.g., \cite{SchJ14,MehU23} for a general overview. 
Recently, a new energy-based formulation has been proposed in \cite{AltS25} which is especially suitable for the modeling of (partial) differential--algebraic equation systems.
This structure arises naturally even in applications which cannot be directly written as a port-Hamiltonian DAE system as introduced in \cite{BeaMXZ18}.
The main motivation of this paper is to investigate the applicability of Krylov subspace methods with short recurrences to linear systems arising from the dissipation-preserving time discretization of linear energy-based DAE systems as introduced in \cite{AltS25,AltCPSS26}.
The main contributions of this paper are listed in the following.
\begin{itemize}[itemsep=0.2em]
	\item We consider the time-discrete system arising from the application of the implicit midpoint rule to a linear energy-based DAE system as introduced in \cite{AltCPSS26} and reformulate it as a linear system with positive semi-definite $\Hmat$ in \Cref{sec:linSystems:modified}.
	\item In \Cref{sec:linSystems:decoupling}, we demonstrate that in the case of a block-diagonal dissipation matrix, also the symmetric part of the resulting linear system is block-diagonal. This further reduces the computational effort for the linear system solves. 
	\item We present numerical results for the biharmonic heat equation and linear poroelasticity in \Cref{sec:biharmonic,sec:poroelasticity}.
	The numerical results for the biharmonic heat equation show that the methods of Widlund and Rapoport typically need significantly less iterations than the generalized minimal residual method (GMRES) and preconditioned GMRES with preconditioneer $\Hmat$.\\	
	For poroelasticity, we consider different choices for solving the linear system with the block-diagonal $\Hmat$ matrix in every iteration and compare the performance with state-of-the-art solvers.
	We observe that the combination of Widlund with a Cholesky decomposition of $\Hmat$ yields the best performance in terms of accuracy and computation time and especially outperforms the preconditioned GMRES method.
	\item We provide a detailed derivation of the method of Widlund (\Cref{alg:widlund}) and a $QR$-based implementation of the method of Rapoport (\Cref{alg:rapoport}) in \Cref{app:Implementation}.
\end{itemize}

The remaining parts of the paper are structured as follows. 
In \Cref{sec:wid_rap} we provide the basic ideas, algorithms, and some theoretical properties of the methods of Widlund and Rapoport.
\Cref{sec:linSystems} is dedicated to the time discretization of energy-based models and the reformulation of the time-discrete system as a linear equation system with positive (semi-)definite symmetric part.
The theoretical findings are then illustrated by means of (numerical) examples in \Cref{sec:decoupling_examples}. Finally, a summary and an outlook are provided in \Cref{sec:conclusion}.
%
%
\section{The methods of Widlund and Rapoport} 
\label{sec:wid_rap}
In this section, we recall the mathematical properties of the methods of Widlund~\cite{Wid78,ConG76} and Rapoport~\cite{Rap78} for solving linear systems of the form $\Amat x = b$ where the symmetric part of the matrix $\Amat$ is positive definite. In order to do so, we follow the more modern description of these methods presented in~\cite{GueLMS22}.

Consider a matrix $\Amat = \Hmat + \Smat\in\R^{n,n}$, where  $\Hmat$ and $\Smat$ are the symmetric and skew-symmetric parts of $\Amat$, respectively, i.e.,
\[
	\Hmat 
	= \tfrac12\, \big( \Amat + \Amat^T \big), \qquad
	\Smat
	= \tfrac12\, \big( \Amat - \Amat^T \big).
\]
Moreover, we assume $\Hmat$ to be positive definite. 

The key idea is to replace the system $\Amat x = b$ with the equivalent system 		
\begin{equation}
\label{eq:prec_sys}
	\big( \Imat_n + \Kmat \big) x 
	= \widehat{b}, \qquad
	\text{where }\
	\Kmat \coloneqq \Hmat^{-1}\Smat
	\ \text{ and } \
	\widehat{b} \coloneqq \Hmat^{-1}b.
\end{equation}
As discussed in \cite{GueLMS22}, the matrix $\Kmat$ is skew symmetric in the inner product induced by the symmetric positive definite matrix $\Hmat$, since  
\begin{equation*}
		\Hmat \Kmat 
		= \Hmat\Hmat^{-1}\Smat 
		= \Smat 
		= \Smat\Hmat^{-1}\Hmat 
		= -\Smat^T\Hmat^{-1}\Hmat 
		= -\Kmat^T \Hmat.
	\end{equation*}
%
Therefore, $\Kmat$ has only purely imaginary eigenvalues and system~\eqref{eq:prec_sys} can be solved via methods based on three-term recurrences. 
For details on the existence of short recurrence methods we refer to \cite[Ch.~4]{LieS12}.
%
%
\subsection{Widlund's method}
\label{sec:wid_rap:wid}

A first method, which can be realized using a three-term recurrence, is the method of Widlund \cite{Wid78} (also derived in \cite{ConG76}). Here, the iterates $x_k$ are defined via 		
\begin{equation}\label{eq:Widlund}
	x_k\in x_0+\calK_k(\Kmat,\widehat{r}_0) 
	\quad \text{such that} \quad 
	r_k	= b - \Amat x_k 
	\perp \calK_k(\Kmat,\widehat{r}_0),
\end{equation}
where $\widehat{r}_0 = \Hmat^{-1}r_0 = \widehat b-(\Imat_n+\Kmat)x_0$ for some starting vector $x_0\in \R^{n}$. Moreover, $\calK_k(\Kmat,\widehat{r}_0)$ denotes the Krylov subspace
\[
	\calK_k(\Kmat,\widehat{r}_0)
	\coloneqq \sspan \big\{ \widehat{r}_0, \Kmat\, \widehat{r}_0, \dots, \Kmat^{k-1}\, \widehat{r}_0 \big\}.
\]
\begin{remark}
Given a general linear system $\Amat x = b$, one usually selects $x_0 = 0$, which yields $r_0 = b$. 
Within this paper, we are especially interested in linear systems arising from a time stepping procedure applied to a linear DAE. This leads to systems of the form $\Amat x^{n+1} = b(x^n)$, where $x^n$ denotes the approximation of the solution at some time point $t^n$. In this case, it is reasonable to set $x^{n}$ as initial guess for the computation of $x^{n+1}$. 
\end{remark}
An implementation of Widlund's method is given in \Cref{alg:widlund}. For completeness, we have included a derivation of the algorithm in \Cref{app:Implementation}.
It was shown in \cite[Thm.~2.1~and~2.2]{Eis83} that the iterates satisfy the optimality properties  
\begin{equation}\label{eq:wid_opt}		
\begin{aligned}
	\norm{x-x_{2k}}_{\Hmat} &= \min_{z\in x_0+(\Imat_n-\Kmat)\calK_{2k}(\Kmat,\widehat{r}_0)}\norm{x-z}_{\Hmat},
	\\ 
	\norm{x-x_{2k+1}}_{\Hmat} &= \min_{z\in x_1+(\Imat_n-\Kmat)\calK_{2k+1}(\Kmat,\widehat{r}_0)}\norm{x-z}_{\Hmat}.
\end{aligned}
\end{equation}
Further, an error bound similar to the one of the CG method can be obtained, as 	
\begin{equation} \label{eq:wid_bound}
	\frac{\norm{x-x_{2k}}_{\Hmat}}{\norm{x-x_0}_{\Hmat}}
	\le 2 \left(\frac{\sqrt{1+\lambda^2}-1}{\sqrt{1+\lambda^2}+1}\right)^k
	\quad \text{and} \quad 
	\frac{\norm{x-x_{2k+1}}_{\Hmat}}{\norm{x-x_1}_{\Hmat}}
	\le 2 \left(\frac{\sqrt{1+\lambda^2}-1}{\sqrt{1+\lambda^2}+1}\right)^k.
\end{equation} 
Here, $\lambda > 0$ is such that the spectrum of the matrix $\Kmat$ is contained in the set  $i\,[-\lambda,\lambda]$, see \cite[Thm.~4.2.]{SzyW93} and \cite[Eq.~5.3]{GueLMS22}.  
These bounds suggest, that smaller values of $\lambda$ should lead to faster convergence of the method. 
If the convergence is slow, a possibility to obtain a faster convergence rate is to use a deflation method as proposed in~\cite{DuFW25}, based on computing approximate eigenspaces corresponding to the eigenvalues of $\Kmat$ with largest magnitude.

\RestyleAlgo{ruled}
\SetKwInOut{Input}{Input}
\SetKwInOut{Output}{Output}
\begin{algorithm}
	\DontPrintSemicolon
	\caption{Widlund's method for solving $\Amat x = b$}
	\label{alg:widlund}
	\Input{$\Amat\in\R^{n,n}$ with $\Hmat\vcentcolon=\tfrac12(\Amat+\Amat^T)$ positive definite,\\ start vector $x_0\in\R^{n}$, right-hand side $b\in\R^{n}$}
	\Output{Approximate solution $x_k$}
	\vspace{0.5em}
	Set $x_{-1}=0$ \; 
	\For{$k=1,2,\dots$}{
		Solve $\Hmat v_k = r_{k-1} = b - \Amat x_{k-1}$ for $v_k$\;
		Set $\rho_k = v_k^T\Hmat v_k = v_k^T r_{k-1}$ \;
		Set $\omega_{k} =
		\begin{cases}
			1, & \text{if }k = 1 \\ 
			(1+ \rho_k/(\rho_{k-1}\omega_{k-1}))^{-1}, & \text{if }k\ge2
		\end{cases} 
		$ \;
		Set $x_k = x_{k-2} +\omega_k(x_{k-1}-x_{k-2}+v_k)$ \;
	}
\end{algorithm}
%
%
\subsection{Rapoport's method}
\label{sec:wid_rap:rap}

A second method for solving the linear system~\eqref{eq:prec_sys} based on a three-term recurrence was described by Rapoport in his PhD thesis \cite{Rap78}. Here, the iterates $x_k$ are defined by 
\begin{equation}\label{eq:Rapoport}
	x_k\in x_0+\calK_k(\Kmat,\widehat{r}_0) 
	\quad \text{such that} \quad 
	r_k = b - \Amat x_k \perp (I_n+\Kmat)\calK_k(\Kmat,\widehat{r}_0).
\end{equation}
An equivalent characterization is given by the minimal residual property 	
\begin{equation}\label{eq:Rapoport_min}
	\norm{r_k}_{\Hmat^{-1}} = \min_{z\in x_0+\calK_k(\Kmat,\widehat{r}_0)}\norm{b-\Amat z}_{\Hmat^{-1}}.
\end{equation}
An upper bound on the $\Hmat^{-1}$-norm of the relative residual is given by 		
\begin{equation} \label{eq:rap_bound}
	\frac{\norm{r_k}_{\Hmat^{-1}}}{\norm{b}_{\Hmat^{-1}}} 
	\le 2 \left(\frac{\lambda}{\sqrt{1+\lambda^2}+1}\right)^{k},
\end{equation}
where, similar as in the method of Widlund, $\lambda>0$ is such that $\sigma(\Kmat)\subseteq i\, [-\lambda,\lambda]$. As for Widlund's method, this bound suggests that we can expect fast convergence for small values of $\lambda$.
A sharper residual bound is presented in \cite[Thm.~2]{DiaFK23}.

A possible implementation of the method of Rapoport is given in \Cref{alg:rapoport} and a derivation of the algorithm is provided in \Cref{app:Implementation}. 
For a given non-zero vector $w$, we use the common short-hand notation $\alpha v = w$ for $\alpha = \norm{w}$ and $v = \tfrac1\alpha w$.

\RestyleAlgo{ruled}
\SetKwInOut{Input}{Input}
\SetKwInOut{Output}{Output}
\begin{algorithm}
    \DontPrintSemicolon
    \caption{Rapoport's method for solving $\Amat x = b$}
    \label{alg:rapoport}
    \Input{$\Amat\in\R^{n,n}$ with $\Hmat\vcentcolon=\tfrac12(\Amat+\Amat^T)$ positive definite,\\ start vector $x_0\in\R^{n}$, right-hand side $b\in\R^{n}$}
    \Output{Approximate solution $x_k$}
	\vspace{0.5em}

    Set $\Kmat=\Hmat^{-1}\Smat$, $\widehat{r}_0 = \Hmat^{-1}b - (\Imat_n+\Kmat)x_0$, 
	$\alpha_0v_1 = \widehat{r}_0$, 
	$\delta_{0} = \alpha_0$, 
    $c_0 = c_{-1} = \gamma_0 = 1$, $s_0=s_{-1}=0$ and
    $v_0 = p_0 = p_{-1}= 0$

    \For{$k=1,2,\dots$}{

    	$\alpha_kv_{k+1} = \Kmat v_k +  \alpha_{k-1} v_{k-1}$

    	$\gamma_{k} = \sqrt{(\gamma_{k-1}c_{k-2})^{2}+\alpha_{k}^2}$\;

		$c_{k} = \gamma_{k-1}c_{k-2}/\gamma_{k}$\;
		
		$s_{k} = \alpha_{k}/\gamma_{k}$\;

		$\delta_k = -s_k\delta_{k-1}$\;

    	$p_{k} =  \left(v_{k} +  \alpha_{k-1}s_{k-2} p_{k-2}\right)/\gamma_{k}$\;

    	$x_{k}  = x_{k-1} + c_{k}\delta_{k-1} p_{k}$ \;
    	
    }
\end{algorithm}
%
%
\section{Linear Systems Arising from Time Discretization}
\label{sec:linSystems}
We focus on the time discretization of energy-based models as introduced in~\cite{AltS25}. In order to guarantee energy dissipation also in the discretized setting, we apply the (implicit) midpoint rule. Moreover, we consider a uniform partition of the time interval with constant step size~$\tau = T/N$, i.e., 
\[
	0 
	= t^0
	< t^1
	< \dots
	< t^N
	= T, \qquad
	t^j = j\tau.
\]
One special case are classical port-Hamiltonian systems, which we discuss first. 
%
%
\subsection{Port-Hamiltonian systems}
\label{sec:linSystems:pH}
Starting point is a system of the form
\begin{align}
	\label{eq:pHDAE}
	\Emat \dot{\state} 
	= (\Jmat - \Rmat) \state + \Bmat \inputPH
\end{align}
with interconnection matrix $\Jmat^T=-\Jmat$, dissipation matrix $\Rmat^T=\Rmat\ge0$, energy matrix $\Emat^T=\Emat\ge0$, and corresponding energy $\hamiltonian(\state) = \frac12 \state^T\Emat\state$. At this point, we would like to emphasize that $\Emat$ may be singular, which would mean that~\eqref{eq:pHDAE} is a DAE. As outlined in~\cite{GueLMS22}, applying the midpoint rule yields the linear system
\[
	\Big[ \Emat - \tfrac\tau2\, (\Jmat - \Rmat) \Big]\, \state^{n+1}
	= \Emat \state^n + \tfrac\tau2\, (\Jmat - \Rmat) \state^n + \tau\, \Bmat \inputPH^{n+1/2}.
\]
The system matrix on the left-hand side decomposes naturally into the symmetric and skew-symmetric parts 
\[
	\Hmat
	= \Emat + \tfrac\tau2 \Rmat, \qquad
	\Smat 
	= - \tfrac\tau2\, \Jmat
\]
with $\Hmat$ being positive semi-definite. 
The special case where $\Hmat$ is even positive definite is of particular relevance and holds, e.g., for ordinary differential equation systems where $\Emat$ is invertible.
In this case, the methods from \Cref{sec:wid_rap} can be directly applied.
Moreover, the skew-symmetric part scales with the step size~$\tau$ such that its influence decreases for refined time grids. 

Further, in \cite[Sec.~4]{GueLMS22} the case of a singular $\Hmat$ is explored. 
Here, the matrix $\Amat=\Hmat+\Smat$ can be transformed to a block diagonal form where the first diagonal block has a positive definite Hermitian part, the last block is skew-symmetric, and the remaining blocks are positive definite.
The transformed system may then be solved efficiently by solving the individual blocks independently, e.g., using one of the methods from \Cref{sec:wid_rap} for the first block, a Krylov method for skew-symmetric problems for the last block (see, e.g.,~\cite{GrePTV16}) and the CG method for the other blocks.
Also direct solvers may be an option, especially for blocks of small or moderate dimension. 
It should be emphasized that the numerical computation of the transformation is challenging for large-scale systems, since it requires several dependent rank-revealing factorizations which is not only costly but also sensitive w.r.t.~perturbations. 
Nevertheless, in practice, the linear system often exhibits a particular structure which may be exploited to derive a transformation to block diagonal form more efficiently or even analytically.
Moreover, the results from \cite[Sec.~4]{GueLMS22} hold for general matrices $\Amat=\Hmat+\Smat$ with symmetric positive semi-definite $\Hmat$ and skew-symmetric $\Smat$ and, hence, may be also applied to the linear systems considered in this paper. 
%
%
\subsection{Energy-based models}
\label{sec:linSystems:energyBased}

We now turn to more general energy-based systems covering an enriched application class. This especially enlarges the classes of constrained systems. Following~\cite{AltS25}, we consider linear systems with quadratic Hamiltonians of the form
\begin{equation*}
	\hamiltonian(\state_1,\state_2) 
	= \tfrac12\, \langle \state_1, \Qmat_1 \state_1 \rangle 
	+ \tfrac12\, \langle \state_2, \Qmat_2 \state_2 \rangle
\end{equation*}
with symmetric positive semi-definite matrices $\Qmat_1,\Qmat_2$. Moreover, we allow a third state variable $\state_3$ which is not part of the energy, leading to systems of the form 
\begin{align*}
	\begin{bmatrix}
		\Qmat_1 \state_1 \\ 
		\dot{\state}_2 \\ 
		0
	\end{bmatrix}
	&=
	(\Jmat - \Rmat)
	\begin{bmatrix}
		\dot{\state}_1\\
		\Qmat_2 \state_2\\
		\state_3
	\end{bmatrix} 
	+ 
	\Bmat \inputPH. 
\end{align*}
It was shown in~\cite{AltS25} that the midpoint rule has the same beneficial properties as for port-Hamiltonian systems, including a power balance on the time-discrete level. 
With the block structure 
\[
	\Jmat
	= \begin{bmatrix}
		\Jmat_{11} & \Jmat_{12} & \Jmat_{13} \\
		\Jmat_{21} & \Jmat_{22} & \Jmat_{23} \\
		\Jmat_{31} & \Jmat_{32} & \Jmat_{33}
	\end{bmatrix}, \qquad
	\Rmat
	= \begin{bmatrix}
		\Rmat_{11} & \Rmat_{12} & \Rmat_{13}\\
		\Rmat_{21} & \Rmat_{22} & \Rmat_{23}\\
		\Rmat_{31} & \Rmat_{32} & \Rmat_{33}
	\end{bmatrix}
\]
the time-discrete system (without inputs) reads
\begin{equation*}
	\begin{bmatrix}
		\Qmat_1\tfrac{\state_1^{n+1}+\state_1^n}2\\
		\tfrac{\state_2^{n+1}-\state_2^n}{\tau}\\
		0
	\end{bmatrix}
	=
	\begin{bmatrix}
		\Jmat_{11}-\Rmat_{11} & \Jmat_{12}-\Rmat_{12} & \Jmat_{13}-\Rmat_{13}\\
		\Jmat_{21}-\Rmat_{21} & \Jmat_{22}-\Rmat_{22} & \Jmat_{23}-\Rmat_{23}\\
		\Jmat_{31}-\Rmat_{31} & \Jmat_{32}-\Rmat_{32} & \Jmat_{33}-\Rmat_{33}
	\end{bmatrix}
	\begin{bmatrix}
		\tfrac{\state_1^{n+1}-\state_1^n}{\tau}\\
		\Qmat_2\tfrac{\state_2^{n+1}+\state_2^n}2\\
		\tfrac{\state_3^{n+1}+\state_3^n}2
	\end{bmatrix}.
\end{equation*}
The resulting linear system to be solved in every time step is hence of the form
\begin{equation*}
	\begin{bmatrix}
		\tfrac{\tau}2 \Qmat_1-(\Jmat_{11}-\Rmat_{11}) & -\tfrac{\tau}2(\Jmat_{12}-\Rmat_{12})\Qmat_2 & -\tfrac{\tau}2(\Jmat_{13}-\Rmat_{13})\\
		-(\Jmat_{21}-\Rmat_{21}) & \Imat_{n_2}-\tfrac{\tau}2(\Jmat_{22}-\Rmat_{22})\Qmat_2 & -\tfrac{\tau}2(\Jmat_{23}-\Rmat_{23})\\
		-(\Jmat_{31}-\Rmat_{31}) & -\tfrac{\tau}2(\Jmat_{32}-\Rmat_{32})\Qmat_2 & -\tfrac{\tau}2(\Jmat_{33}-\Rmat_{33})
	\end{bmatrix}
	\begin{bmatrix}
		\state_1^{n+1}\\
		\state_2^{n+1}\\
		\state_3^{n+1}
	\end{bmatrix}
	= b(\state_1^n,\state_2^n,\state_3^n).
\end{equation*}
In contrast to the classical port-Hamiltonian systems considered before, the symmetric part of the system matrix is not guaranteed to be positive semi-definite.
One reason appears to be $\Qmat_2$ occurring in the second block column and destroying the structure.  
%
%
\subsection{Modification of the system}
\label{sec:linSystems:modified}
To circumvent the observed issue, inspired by \cite{GueLMS22} (among others), we replace the original system by 
\begin{align}
	\label{eq:energyBasedModel:withE2}
	\begin{bmatrix}
		\Qmat_1 \state_1 \\ 
		\Emat_2 \dot{\state}_2 \\ 
		0
	\end{bmatrix}
	=
	(\Jmat - \Rmat)
	\begin{bmatrix}
		\dot{\state}_1\\
		\state_2\\
		\state_3
	\end{bmatrix} 
	+ 
	\Bmat \inputPH 
\end{align}
with associated Hamiltonian
\begin{equation*}
	\hamiltonian(\state_1,\state_2) 
	= \tfrac12\, \langle \state_1, \Qmat_1 \state_1 \rangle 
	+ \tfrac12\, \langle \state_2, \Emat_2 \state_2 \rangle,
\end{equation*}
where $\Qmat_1,\Emat_2$ are again symmetric positive semi-definite. 
Note that this is a special case of the extension presented in~\cite{AltCPSS26}. Hence, also for this system the midpoint rule guarantees a discrete power balance. The resulting time-discrete system reads
\begin{equation*}
	\begin{bmatrix}
		\Qmat_1\tfrac{\state_1^{n+1}+\state_1^n}2\\
		\Emat_2\tfrac{\state_2^{n+1}-\state_2^n}{\tau}\\
		0
	\end{bmatrix}
	=
	\begin{bmatrix}
		\Jmat_{11}-\Rmat_{11} & \Jmat_{12}-\Rmat_{12} & \Jmat_{13}-\Rmat_{13}\\
		\Jmat_{21}-\Rmat_{21} & \Jmat_{22}-\Rmat_{22} & \Jmat_{23}-\Rmat_{23}\\
		\Jmat_{31}-\Rmat_{31} & \Jmat_{32}-\Rmat_{32} & \Jmat_{33}-\Rmat_{33}
	\end{bmatrix}
	\begin{bmatrix}
		\tfrac{\state_1^{n+1}-\state_1^n}{\tau}\\
		\tfrac{\state_2^{n+1}+\state_2^n}2\\
		\tfrac{\state_3^{n+1}+\state_3^n}2
	\end{bmatrix},
\end{equation*}
where we use the same block structure of $\Jmat$ and $\Rmat$ as before. Multiplying everything by $\tau$, we can write the linear system as 
\begin{align*}
	\Bigg(
	\begin{bmatrix}
		\tfrac{\tau}2 \Qmat_1 & \Nullmat & \Nullmat\\
		\Nullmat & \tfrac2{\tau} \Emat_2 & \Nullmat\\
		\Nullmat & \Nullmat & \Nullmat
	\end{bmatrix}
	- \Jmat + \Rmat
	\Bigg)
	\begin{bmatrix}
		\state_1^{n+1}\\
		\tfrac{\tau}2\state_2^{n+1}\\
		\tfrac{\tau}2\state_3^{n+1}
	\end{bmatrix}
	= b(\state_1^n,\state_2^n,\state_3^n).
\end{align*}
Therefore, when considering $[\state_1^{n+1};\tfrac{\tau}2\state_2^{n+1};\tfrac{\tau}2\state_3^{n+1}]$ as the vector of unknowns, the system matrix of the linear system is the sum of a symmetric positive semi-definite and a skew-symmetric matrix, similarly as in \cite{GueLMS22}.
%
%
\subsection{Decoupling of the $\Hmat$ matrix}
\label{sec:linSystems:decoupling}
As seen in Section~\ref{sec:wid_rap:wid}, Widlund's three-term iteration only requires the solution of linear systems with the symmetric and positive definite matrix~$\Hmat$. The same holds true for the method of Rapoport from Section~\ref{sec:wid_rap:rap}. We are interested in applications where this matrix is block diagonal such that the solution of $\Hmat x = b$ decouples into a number of smaller subsystems, all of them again with a symmetric positive definite matrix. This happens if the dissipation matrix $\Rmat$ is block diagonal, leading to 
\[
	\Hmat 
	= 
	\begin{bmatrix}
		\tfrac{\tau}2 \Qmat_1 + \Rmat_{11} & \Nullmat & \Nullmat\\
		\Nullmat & \tfrac2{\tau} \Emat_2 + \Rmat_{22} & \Nullmat\\
		\Nullmat & \Nullmat & \Rmat_{33}
	\end{bmatrix}.
\]
In this case, the dissipation matrix is decoupled and, therefore, the coupling of $\state_1$, $\state_2$, and~$\state_3$ in~\eqref{eq:energyBasedModel:withE2} is only due to the off-diagonal blocks in the interconnection matrix $\Jmat$.
In the remaining parts of the paper, we illustrate this decoupling for several applications. 
%
%
\section{Decoupling Examples}
\label{sec:decoupling_examples}
Within this section, we collect (constrained) examples which fit into the energy-based framework of Section~\ref{sec:linSystems:modified} and decouple as outlined in Section~\ref{sec:linSystems:decoupling}. Classical port-Hamiltonian (DAE) examples 
are already covered by the theory presented in~\cite{GueLMS22} and are not recapitulated in detail here. 
As preparation for the biharmonic heat equation, we start with an application without constraint. Afterwards, we discuss the equations of linear poroelasticity.
%
%
\subsection{Beam equation} 
We consider the fourth-order equation $\ddot u + \Delta^2 u = f$, which is used, e.g., for describing the dynamics of an Euler--Bernoulli beam, see e.g.~\cite[Sec.~6.5]{Inm14}. 
To obtain a first-order formulation of the system (w.r.t.~time) we introduce $v\vcentcolon= \dot{u}$ and $w \vcentcolon= \Delta u$. This leads to the system  
\begin{align*}
	\dot v + \Delta w 
	&= f, \\
	\dot{w} - \Delta v 
	&= 0.	
\end{align*}
Note that the original solution $u$ is not part of the system but can be recovered, e.g., by integrating $v$ over time. Considering homogeneous Dirichlet boundary conditions and a standard finite element discretization, this then leads to the two equations
\begin{align*}
	\Mmat\dot v - \Kmat w 
	&= f, \\
	\Mmat \dot{w} + \Kmat v 
	&= 0
\end{align*}
with a mass matrix $\Mmat$ and a stiffness matrix $\Kmat$. Note that both matrices are symmetric and positive definite. The corresponding total energy is the sum of the kinetic and the potential energy, given by 
\[
	\hamiltonian
	= \tfrac12\, v^T \Mmat v + \tfrac12\, w^T \Mmat w.
\]
Hence, considering $\state_2 = \begin{bsmallmatrix}	v\\ w \end{bsmallmatrix}$ in the notion of Section~\ref{sec:linSystems:modified}, we have $\hamiltonian = \frac12 \state_2^T\Emat_2 \state_2$ with $\Emat_2 = \diag(\Mmat, \Mmat)$. 
The corresponding system reads 
\[
	\Emat_2\dot \state_2
	= 
	\begin{bmatrix}
		\Mmat & \Nullmat \\
		\Nullmat & \Mmat
	\end{bmatrix}
	\begin{bmatrix}
		\dot v\\
		\dot w
	\end{bmatrix}
	=
	\begin{bmatrix} 
		\phantom{-}\Nullmat & \Kmat \\ 
		-\Kmat & \Nullmat
	\end{bmatrix}
	\begin{bmatrix} v \\ w \end{bmatrix}
	+
	\begin{bmatrix} f \\ 0 \end{bmatrix}
	=
	\Jmat\state_2+\begin{bmatrix} f \\ 0 \end{bmatrix},
\]
which is of the form \eqref{eq:energyBasedModel:withE2} with empty $\state_1$ and $\state_3$. In particular, it has the port-Hamiltonian structure considered in~\cite{GueLMS22}. 
An application of the midpoint rule then yields a system with iteration matrix $\Hmat+\Smat$ with 
%
%
\[
	\Hmat 
	= \frac2\tau \begin{bmatrix} 
		\Mmat & \Nullmat\\ 
		\Nullmat & \Mmat 
	\end{bmatrix}, \qquad 
	\Smat
	= \begin{bmatrix} 
		\Nullmat & -\Kmat\\ 
		\Kmat & \phantom{-}\Nullmat 
	\end{bmatrix}.
\]
Note that, in this example, $\Hmat$ is only a suitable preconditioner if the time step size $\tau$ is sufficiently small. For moderate $\tau$, on the other hand, the skew-symmetric part dominates as it contains the stiffness matrix.  
%
%
\subsection{Biharmonic heat equation} 
\label{sec:biharmonic}
\newcommand{\nnu}{\eta}
The second example is devoted to the fourth-order linear diffusion equation $\dot u + \Delta^2 u = f$, also called the biharmonic heat equation. This equation is also connected to the (linear part) of the extended Fisher--Kolmogorov equation~\cite{PelT97}, where the fourth-order term serves as stabilization term within the pattern formation in bi-stable systems \cite{DeeS88}. As coupled system of second order, we can write this as  
\begin{equation}\label{eq:biharmonic_heat_eq}
\begin{aligned}
	\dot u - \Delta w 
	&= f, \\
	\Delta u + w 
	&= 0.	
\end{aligned}
\end{equation}
Considering homogeneous Dirichlet boundary conditions for $u$ and $w$ and a finite element discretization as before, this leads to the semi-discrete system 
\begin{align*}
	\Mmat \dot u + \Kmat w
	&= f, \\
	-\Kmat u + \Mmat w 
	&= 0 
\end{align*}
with mass matrix $\Mmat\in\R^{\nnu,\nnu}$ and stiffness matrix $\Kmat\in\R^{\nnu,\nnu}$. Again, both matrices are symmetric positive definite. In the one-dimensional setting with spatial domain $(0,1)$ and first-order finite elements, we have $\nnu+1 = \frac 1h$, where $h$ denotes the considered mesh size. We present two ways how this system fits into the energy-based framework of Section~\ref{sec:linSystems:modified}.  
\medskip

First, we consider $\state_2 = u$, $\state_3 = w$ together with the energy function $\hamiltonian = \frac12 u^T\Emat_2 u$ and $\Emat_2 = \Mmat$. This leads to the system 
\[
	\begin{bmatrix} \Mmat\dot u \\ 0	\end{bmatrix}
	= 
	\begin{bmatrix} \Nullmat & -\Kmat \\ \Kmat & -\Mmat \end{bmatrix}
	\begin{bmatrix} u \\ w \end{bmatrix}
	+
	\begin{bmatrix} f \\ 0 \end{bmatrix}
	=
	\Bigg(
	\begin{bmatrix} \Nullmat & -\Kmat \\ \Kmat & \phantom{-}\Nullmat \end{bmatrix}
	-
	\begin{bmatrix} \Nullmat & \Nullmat \\ \Nullmat & \Mmat \end{bmatrix}
	\Bigg)
	\begin{bmatrix} u \\ w \end{bmatrix}
	+
	\begin{bmatrix} f \\ 0 \end{bmatrix},
\]
which is of the form \eqref{eq:energyBasedModel:withE2} with empty $z_1$. 
Within this setting, the application of the midpoint rule yields the symmetric and skew-symmetric parts  
\begin{equation}\label{eq:biharmonic_v1}
	\Hmat_1 
	= \begin{bmatrix} 
		\frac{2}{\tau}\Mmat & \Nullmat\\ 
		\Nullmat & \Mmat 
	\end{bmatrix}, \qquad 
	\Smat_1 
	= \begin{bmatrix} 
		\phantom{-}\Nullmat & \Kmat \\
		-\Kmat & \Nullmat 
	\end{bmatrix}.
\end{equation}
Having the properties of the finite element matrices in mind, one can already see that for moderate step sizes (or balanced discretizations in time and space) the symmetric part is not dominant in this formulation.
\medskip

Second, we consider $\state_1 = u$, $\state_3 = w$, and the energy function $\hamiltonian = \frac12 u^T\Qmat_1 u$ with $\Qmat_1 = \Kmat$. In this case, we obtain the system 
\[
	\begin{bmatrix} \Kmat u \\ 0	\end{bmatrix}
	= 
	\begin{bmatrix} \phantom{-}\Nullmat & \phantom{-}\Mmat \\ -\Mmat & -\Kmat \end{bmatrix}
	\begin{bmatrix} \dot u \\ w \end{bmatrix}
	+
	\begin{bmatrix} 0 \\ f \end{bmatrix}
	=
	\Bigg(
	\begin{bmatrix} \phantom{-}\Nullmat & \Mmat \\ -\Mmat & \Nullmat \end{bmatrix}
	-
	\begin{bmatrix} \Nullmat & \Nullmat\\ \Nullmat & \Kmat \end{bmatrix}
	\Bigg)
	\begin{bmatrix} \dot u \\ w \end{bmatrix}
	+
	\begin{bmatrix} 0 \\ f \end{bmatrix},
\]
which is again of the form~\eqref{eq:energyBasedModel:withE2}. The midpoint rule (with the rescaling from Section~\ref{sec:linSystems:modified}) yields the system
%
%
%
\begin{equation}\label{eq:biharmonic_v2_sys}
	\begin{bmatrix} 
		\frac{\tau}{2}\Kmat & -\Mmat \\ 
		\Mmat &  \phantom{-}\Kmat 
	\end{bmatrix}
	\begin{bmatrix} u^{n+1} \\  \frac{\tau}{2} w^{n+1} \end{bmatrix} 
	= b(u^n, w^n)
	\coloneqq 
	\begin{bmatrix} 
		- \frac{\tau}{2} \Kmat u^n + \frac{\tau}{2} \Mmat w^n \\ 
		\Mmat u^n - \frac{\tau}{2} \Kmat w^n + \tau f(t^{n+1/2})
	\end{bmatrix}.
\end{equation}
In this formulation, the symmetric and skew-symmetric parts are given by 
\begin{equation}\label{eq:biharmonic_v2}	
	\Hmat_2 
	= \begin{bmatrix} 
		\frac{\tau}{2}\Kmat & \Nullmat\\ 
		\Nullmat & \Kmat 
	\end{bmatrix}, \qquad 
	\Smat_2 
	= \begin{bmatrix} 
		\Nullmat & -\Mmat \\ 
		\Mmat & \phantom{-}\Nullmat 
	\end{bmatrix}.
\end{equation}
In order to compare the two formulations, we have computed the maximum absolute value of the eigenvalues of $\Hmat_j^{-1}\Smat_j$ for $j =1,2$ using the Matlab routine \texttt{eig}. The size of the matrices is usually coupled to the time step size via $\tau \approx h\approx \nnu^{-1}$ or $\tau \approx \sqrt{h}\approx \nnu^{-1/2}$. 
The results for these choices of $\tau$ and $\nnu$ are presented in \Cref{tab:lambda_bound}.
\begin{table}
\begin{tabular}{c||*{3}{c}||*{3}{c}}
	 & \multicolumn{3}{c||}{$\tau = \nnu^{-1}$}& \multicolumn{3}{c}{$\tau = \nnu^{-\frac{1}{2}} $} \\ 
	$\nnu$ & $\tau$ & $\displaystyle\max_{\lambda\in\sigma(\Hmat_1^{-1}\Smat_1)}\abs{\lambda}$ & $\displaystyle\max_{\lambda\in\sigma(\Hmat_2^{-1}\Smat_2)}\abs{\lambda}$ & $\tau$ & $\displaystyle\max_{\lambda\in\sigma(\Hmat_1^{-1}\Smat_1)}\abs{\lambda}$ & $\displaystyle\max_{\lambda\in\sigma(\Hmat_2^{-1}\Smat_2)}\abs{\lambda}$ \\ 
	\hline\hline
	$10^{1}$ & $10^{-1}$ & $3.06 \cdot 10^{2}$ & $4.50 \cdot 10^{-1}$ & $10^{-\frac{1}{2}}$ & $5.44 \cdot 10^{2}$ & $2.53 \cdot 10^{-1}$ \\
	$10^{2}$ & $10^{-2}$ & $8.65 \cdot 10^{3}$ & $1.43 \cdot 10^{0}$ & $10^{-1}$ & $2.74 \cdot 10^{4}$ & $4.53 \cdot 10^{-1}$ \\
	$10^{3}$ & $10^{-3}$ & $2.69 \cdot 10^{5}$ & $4.53 \cdot 10^{0}$ & $10^{-\frac{3}{2}}$ & $1.51 \cdot 10^{6}$ & $8.06 \cdot 10^{-1}$ \\
\end{tabular}
\caption{Largest absolute value of the eigenvalues of $\Hmat_j^{-1}\Smat_j$ for different time step sizes $\tau$ and matrix dimensions $\nnu$.}
\label{tab:lambda_bound}
\end{table}
In view of the convergence bounds~\eqref{eq:wid_bound} and~\eqref{eq:rap_bound}, a smaller value of $\lambda$ is preferable. 
Since the absolute value of the largest eigenvalue of $\Hmat_1^{-1}\Smat_1$ is much larger than the one of $\Hmat_2^{-1}\Smat_2$ for both choices of $\tau$, these bounds indicate that the second formulation with the matrices~\eqref{eq:biharmonic_v2} is favorable.

Next, we compare the performance of the methods of Widlund and Rapoport for solving the linear system \eqref{eq:biharmonic_v2_sys} resulting from the first time step of the second formulation, to the performance of GMRES and preconditioned GMRES when applied to this system. 
Here, we choose $\Hmat_2$ as preconditioner such that the preconditioned GMRES minimizes the original residual in the $\Hmat_2^{-2}$-norm. 
We emphasize that there are more sophisticated preconditioning strategies, but for our comparison, the main point is that for this choice, preconditioned GMRES and the methods of Widlund and Rapoport are all (implicitly) solving the system \eqref{eq:prec_sys}, while (unpreconditioned) GMRES acts as a reference to the original system.  
Here we choose 
\[
	u(0,x) = \sin(\pi x) \text{ for } x\in [0,1], \qquad 
	f(t) = t,
\]
and $\tau=\nnu^{-1}$ for different step sizes $\tau$. Additionally, we choose the initial vector to be the zero vector, since further experiments indicated that there is no substantial advantage in choosing the initial vector as, e.g., the solution of an earlier time step due to the small number of iterations needed.

We solve all systems until the (computed) relative residual norm is less than $10^{-6}$. 
Since the involved matrices become increasingly ill-conditioned with growing $\nnu$, a lower tolerance was not achievable in our experiments.  
In each iteration, the linear system involving $\Hmat_2$ is solved via a precomputed Cholesky factorization. 
Since $\Hmat_2$ is tridiagonal, the full Cholesky factor is given by a bidiagonal matrix and, hence, can be computed without losing sparsity. 
The convergence behavior of the methods is depicted in \Cref{fig:biharmonic_single_sys}. The corresponding computation times of the different solvers, as well as the time it took to compute the Cholesky decomposition of $\Hmat_2$ are reported in \Cref{tab:biharmonic_single_sys}. The timings for GMRES are not shown here, as the method was not able to reduce the residual below 
$10^{-6}$ within $100$ iterations.
For all examined values of $\nnu$, we observe that the methods of Widlund and Rapoport only need a few steps to reduce the norm of the residual below the specified tolerance. The performance of preconditioned GMRES, on the other hand, becomes increasingly slower for larger values of $\nnu$. 
In terms of computation times, we observe that the methods of Widlund and Rapoport perform better than preconditioned GMRES. This has to be expected as the former methods rely on short recurrences, while GMRES uses a full recurrence. 

\begin{figure}
    \centering
    \input{pics/biharmonic_widlund_norms.tex}
    \hfill
    \input{pics/biharmonic_rapoport_norms.tex}
    \hfill
    \input{pics/biharmonic_pgmres_norms.tex}
    \small
\begin{tikzpicture}
    \begin{axis}[
      hide axis,
      scale only axis,
      xmin=0, xmax=1,
      ymin=0, ymax=1,
      width=0pt,
      height=0pt,
      legend columns=4,
      legend style={
        draw=black,
        /tikz/every even column/.append style={column sep=0.5cm}
      }
    ]

    \addlegendimage{very thick, mycolor1}
    \addlegendentry{$\norm{\,\cdot\,}_2$}

    \addlegendimage{very thick, mycolor2, dashed}
    \addlegendentry{$\norm{\,\cdot\,}_{\Hmat^{-1}}$}

    \addlegendimage{very thick, mycolor5, dash dot}
    \addlegendentry{$\norm{\,\cdot\,}_{\Hmat^{-2}}$}

    \end{axis}
    \end{tikzpicture}
    \caption{Comparison of the three different norms of the relative residual when solving system~\eqref{eq:biharmonic_v2_sys} with $\nnu = \tau^{-1} = 10^5$.}
	\label{fig:norms_biharmonic}
\end{figure}
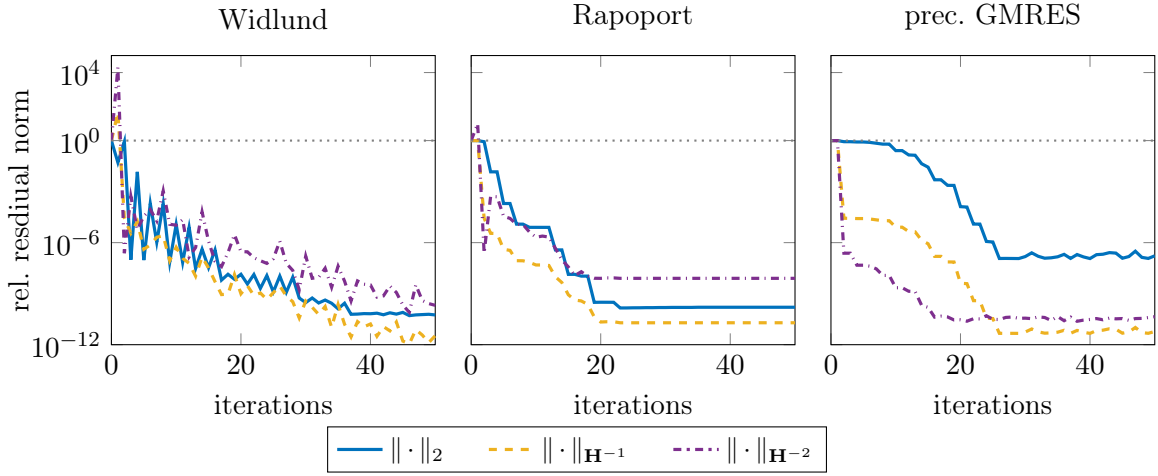

We would like to emphasize that none of the methods is actually minimizing the Euclidean norm. Instead, as discussed previously, the method of Rapoport minimizes the residual in the  $\Hmat^{-1}$-norm and preconditioned GMRES minimizes the residual in the $\Hmat^{-2}$-norm, while the method of Widlund does not minimize the residual norm, but rather the norm of the error. 
In \Cref{fig:norms_biharmonic}, the convergence behaviour in the Euclidean norm, the $\Hmat^{-1}$-norm and the $\Hmat^{-2}$-norm are compared for
system~\eqref{eq:biharmonic_v2_sys} with $\nnu = \tau^{-1}=10^5$ for $50$ iterations. 
To reduce the loss of orthogonality in preconditioned GMRES (especially once the iteration is stagnating), we included a reorthogonalization step. 
We can see that the different norms seem to behave similarly for the methods of Widlund and Rapoport. For preconditioned GMRES, on the other hand, there is a large difference between the Euclidean norm and the other two norms. This also explains the final accuracy of $10^{-6}$ (in the Euclidean norm), as preconditioned GMRES has already reached its final accuracy in the $\Hmat^{-2}$-norm.
\begin{remark}
Further experiments, involving a predefined solution and an accordingly constructed right-hand side, have indicated that even though the methods of Widlund and Rapoport, as well as preconditioned GMRES, reduce the relative residual norm below $10^{-6}$, the final error can behave quite differently for the three methods. 
In particular, the final error can be larger for the methods of Widlund and Rapoport, when compared to preconditioned GMRES. This may be caused by the higher number of iterations performed by preconditioned GMRES. 
\end{remark}
\begin{remark}
In the examples considered here, the linear systems are of unsymmetric saddle point form. In fact, multiplying the second block row by $-1$ yields a symmetric linear system with saddle point structure. 
Both, the symmetric (see, e.g., the survey article~\cite{BenGL05}) as well as the unsymmetric formulation (see, e.g., \cite{BenS06,LieP08,LieR24}) have been studied extensively in the literature.
However, since this structure is a feature of the contemplated examples and not of the model class and its discretization itself, we have not considered methods specifically tailored to (unsymmetric) saddle point systems in this work.
\end{remark}
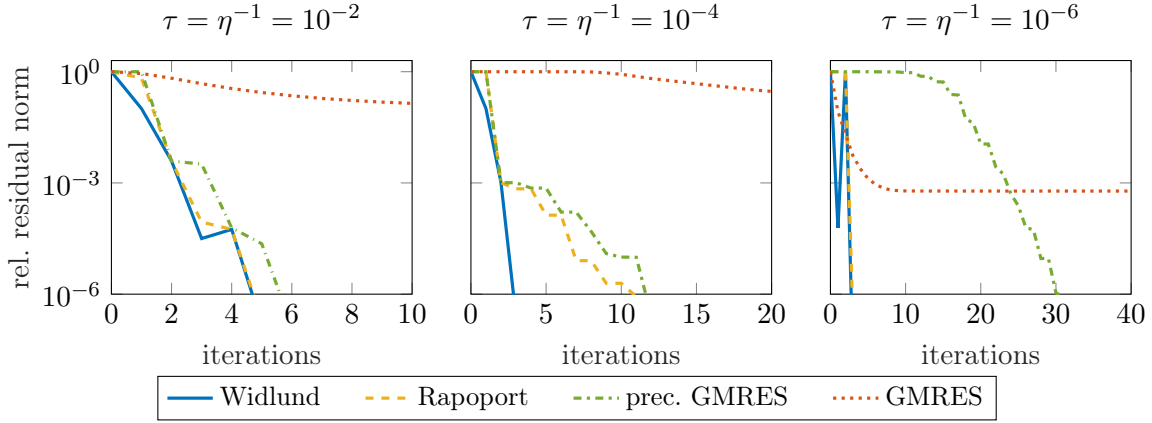
\begin{figure}
    \centering
    \input{pics/biharmonic_1e-2.tex}
    \hfill
    \input{pics/biharmonic_1e-4.tex}
    \hfill
    \input{pics/biharmonic_1e-6.tex}
    \small 
    \begin{tikzpicture}
    \begin{axis}[
      hide axis,
      scale only axis,
      xmin=0, xmax=1,
      ymin=0, ymax=1,
      width=0pt,
      height=0pt,
      legend columns=4,
      legend style={
        draw=black,
        /tikz/every even column/.append style={column sep=0.5cm}
      }
    ]

    \addlegendimage{very thick, mycolor1}
    \addlegendentry{Widlund}

    \addlegendimage{very thick, mycolor2, dashed}
    \addlegendentry{Rapoport}

    \addlegendimage{very thick, mycolor3, dash dot}
    \addlegendentry{prec.~GMRES}

    \addlegendimage{very thick, mycolor4, dotted}
    \addlegendentry{GMRES}

    \end{axis}
    \end{tikzpicture}
    \caption{Relative residual norms for the four methods and different time step sizes $\tau$ and system dimensions $\nnu$.} 
    \label{fig:biharmonic_single_sys}
\end{figure}
\begin{table}
\begin{tabular}{c||c|c|c|c|c}
	$\tau = \nnu^{-1} $ & $10^{-2} $ & $10^{-3} $ & $10^{-4} $ & $10^{-5} $ & $10^{-6} $ \\[0.1em]
	\hline \hline
	Cholesky  & 
	$1.120 \cdot 10^{-4} $ & $1.770 \cdot 10^{-4} $ & $1.014 \cdot 10^{-3} $ & $1.257 \cdot 10^{-2} $ & $1.693 \cdot 10^{-1} $ \\
	\hline 
	\hline
	Widlund  & 
	$2.871 \cdot 10^{-3} $ & $2.227 \cdot 10^{-3} $ & $3.830 \cdot 10^{-3} $ & $4.056 \cdot 10^{-2} $ & $5.087 \cdot 10^{-1} $ \\
	\hline 
	Rapoport  & 
	$4.124 \cdot 10^{-3} $ & $1.836 \cdot 10^{-3} $ & $1.238 \cdot 10^{-2} $ & $1.417 \cdot 10^{-1} $ & $5.455 \cdot 10^{-1} $ \\
	\hline 
	prec.~GMRES  & 
	$5.701 \cdot 10^{-3} $ & $2.654 \cdot 10^{-3} $ & $1.178 \cdot 10^{-2} $ & $3.760 \cdot 10^{-1} $ & $6.403 \cdot 10^{0} $ \\
\end{tabular}
\caption{Computation times for the Cholesky decomposition of $\Hmat_2$ and the solution of system~\eqref{eq:biharmonic_v2_sys} using the methods of Widlund, Rapoport, and preconditioned GMRES for different time step sizes $\tau$ and dimension~$\nnu = \tau^{-1}$.}
\label{tab:biharmonic_single_sys}
\end{table}
%
%
%
\subsection{Linear poroelasticity} 
\label{sec:poroelasticity}
Going back to~\cite{Bio41,Sho00}, the (spatially discretized) equations of linear poroelasticity are given by the DAE 
\begin{align}
\label{eq:poro:semiDiscrete}
	\begin{bmatrix} \Nullmat & \Nullmat \\ \Dmat & \Cmat \end{bmatrix}
	\begin{bmatrix} \dot u \\ \dot p \end{bmatrix}
	= 
	\begin{bmatrix} -\Amat & \phantom{-}\Dmat^T \\ \phantom{-}\Nullmat & -\Bmat\phantom{^T} \end{bmatrix}
	\begin{bmatrix} u \\ p \end{bmatrix}
	+ 
	\begin{bmatrix} \f \\ \g \end{bmatrix}
\end{align}
with initial conditions $u(0) = u^0 \in \R^n$ and $p(0) = p^0 \in \R^m$. 
In applications, $u$ models the deformation of the porous medium whereas $p$ equals the pressure of the incompressible viscous fluid with which the material is saturated. The corresponding energy has the form
\[
	\hamiltonian
	= \tfrac12\, u^T \Amat u + \tfrac12\, p^T \Cmat p.
\] 
In the system equations, $\Amat\in\R^{n,n}$ equals the stiffness matrix from linear elasticity based on the Lam{\'e} coefficients $\lambda$ and $\mu$. The matrix $\Bmat\in\R^{m,m}$ is the diffusion matrix, where the diffusion coefficient equals the permeability~$\kappa$ divided by the fluid viscosity~$\nnu$. Further, $\Cmat\in\R^{m,m}$ is a mass matrix scaled with one over the Biot modulus $M$ and $\Dmat\in\R^{m,n}$ is a rectangular matrix which scales with the Biot--Willis ﬂuid--solid coupling coefficient~$\alpha$. 
In the upcoming experiments, we consider first-order finite elements, which are realized in Matlab and based on~\cite{AlbCFK02}. 
Further, we set 
\[
	\lambda = \mu = 10,\qquad
	\kappa/\nnu = 1,\qquad 
	M = 10,\qquad
	\alpha = 1
\]
for the physical parameters.
%
%
\subsubsection{Midpoint rule}
The application of any implicit time stepping scheme to~\eqref{eq:poro:semiDiscrete} yields a large (coupled) linear system, which has to be solved in every time step. 
Considering the midpoint rule with $\frac\tau2 p$ as pressure variable as proposed in Section~\ref{sec:linSystems:modified}, the resulting time stepping scheme reads 
\[
	\begin{bmatrix} 
		\frac\tau2 \Amat & -\Dmat^T\\ 
		\Dmat & \frac2\tau \Cmat + \Bmat 
	\end{bmatrix}
	\begin{bmatrix} 
		u^{n+1} \\ 
		\frac\tau2 p^{n+1} 
	\end{bmatrix}
	=  
	b(u^n, p^n)
	\coloneqq
	\begin{bmatrix} 
		- \frac\tau2 \Amat u^n + \frac\tau2 \Dmat^T p^n + \tau f^{n+1/2} \\ 
		\Cmat p^n + \Dmat u^n - \frac\tau2 \Bmat p^n + \tau g^{n+1/2} 
	\end{bmatrix}.
\]
We would like to emphasize that this scheme is second-order convergent 
and dissipation-preserving. The iteration matrix can be decomposed into $\Hmat + \Smat$ with
\[
	\Hmat 
	= \begin{bmatrix}
		\tfrac{\tau}2 \Amat & \Nullmat \\
		\Nullmat & \tfrac2{\tau}\Cmat+\Bmat
	\end{bmatrix}, \qquad
	\Smat
	= \begin{bmatrix}
		\Nullmat & -\Dmat^T \\
		\Dmat & \Nullmat
	\end{bmatrix}.
\]
The properties of $\Amat$ and $\Cmat$ imply that $\Hmat$ is positive definite such that Widlund's as well as Rapoport's method are applicable. In the following experiments, however, we only consider the method of Widlund. Since $\Amat$ and $\Bmat$ result from differential operators of second order and $\Dmat$ from a divergence operator, one can see that $\Hmat$ is dominant as long as $\tau$ does not get too small. 

%
%
In the first experiment, we study the computation times for a single time step. As computational domain, we consider the unit square with a uniform triangulation leading to system dimensions $n = 130{,}050$ and $m = 65{,}025$. The time step size is set to $\tau = \sqrt{2}\cdot2^{-8}$, which equals the spatial mesh size. As right-hand sides, we use $f(t) = g(t) = 1$. The initial data is set to $p(0) = \sin(\pi x) \sin(\pi y)$, which is consistent with the assumed homogeneous Dirichlet boundary conditions. 

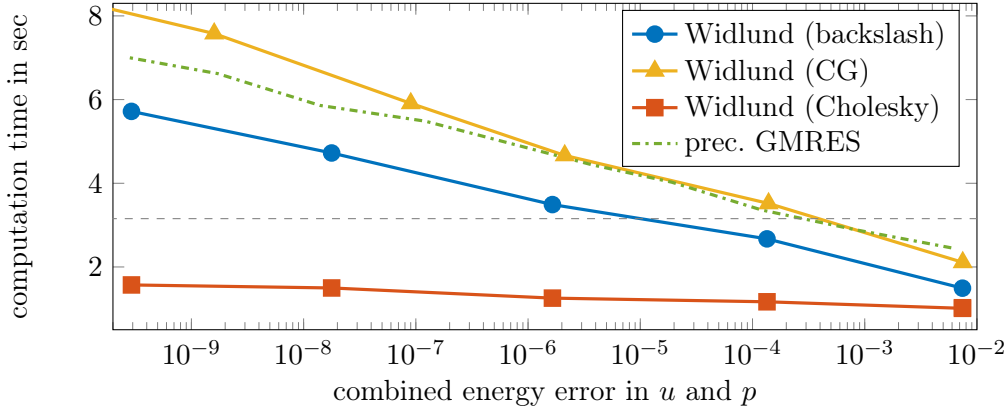
\begin{figure}
	\input{pics/square_oneStep_nRef8.tex} 
	\caption{Computation times vs.~energy error for a single time step of the midpoint rule using different iteration schemes. The gray dashed line indicates the $3.15$ seconds which were needed to compute the reference solution with Matlab's backslash operator applied to the overall system. }
	\label{fig:poro:oneTimeStep}
\end{figure}

The computation time of \Cref{alg:widlund} depends on the implemented stopping criterion. We prescribe different tolerances~$\tol$ for the relative residual and terminate the computation as soon as 
\[
	\frac{\|b(u^n, p^n) - (\Hmat+\Smat) x_k\| }{\|b(u^n, p^n)\|} 
	\le \tol.
\] 
As starting vector, we implement the initial data.  
Recall that in each iteration of Widlund's method, a linear system involving $\Hmat$ has to be solved and, since $\Hmat$ is block diagonal, this system decouples into two smaller subsystems. 
For the solution of these subsystems, we compare the standard backslash operator from Matlab with a CG iteration (with inner tolerance $\frac14 \tol$) as well as a Cholesky decomposition (with minimum degree ordering; see~\cite{AmeDD04}) 
of the matrices~$\Amat$ and $\frac2\tau\Cmat+\Bmat$, respectively. The resulting computation times are presented in Figure~\ref{fig:poro:oneTimeStep}. The indicated energy error equals the error in $u$ measured in the $\Amat$-norm plus the error in $p$ measured in the $\Cmat$-norm. As comparison, we also show the results of a GMRES iteration applied to the full system with $\Hmat$ as a preconditioner.  
%
%
\subsubsection{Comparison to other decoupling methods}
There are many other approaches that aim to decouple the resulting system equations, e.g., the well-known fixed stress iteration~\cite{KimTJ11a}. This scheme is based on the implicit Euler discretization together with a fixed-point iteration. A generalization to higher order based on backward differentiation formulae (BDF) is given in~\cite{AltMU24b}. These schemes require a stabilization parameter~$L$ and a maximum number of inner iteration steps~$K$. It was shown that convergence is only guaranteed if $K$ grows for decreasing time step sizes~$\tau$. A similar approach was introduced in~\cite{AltD25}. Here, however, the number of inner iteration steps~$K$ only depends on the physical parameters, leading to a scheme that converges with order two. This provides a very efficient scheme but does not guarantee a preservation of the energy dissipation. 
Finally, we include the method from~\cite{AltMU24} into the comparison, a method without any inner iteration, which only converges in specific parameter settings. 
%
%
%

\begin{figure}
	\input{pics/square_times_nRef6_T1_om05} 
	\caption{Comparison of different decoupling time stepping methods of order two in terms of computation time. The fixed stress scheme is shown for $L=0$ and $K=3$ (dotted), $K=6$ (dashed), and $K=9$ (solid). Widlund's method is run with $\tol = \tau^2$ (dotted), $\tol = 0.01\,\tau^2$ (dashed), and $\tol = 0.001\,\tau^2$ (solid).}
	\label{fig:poro:computationTimes}
\end{figure}
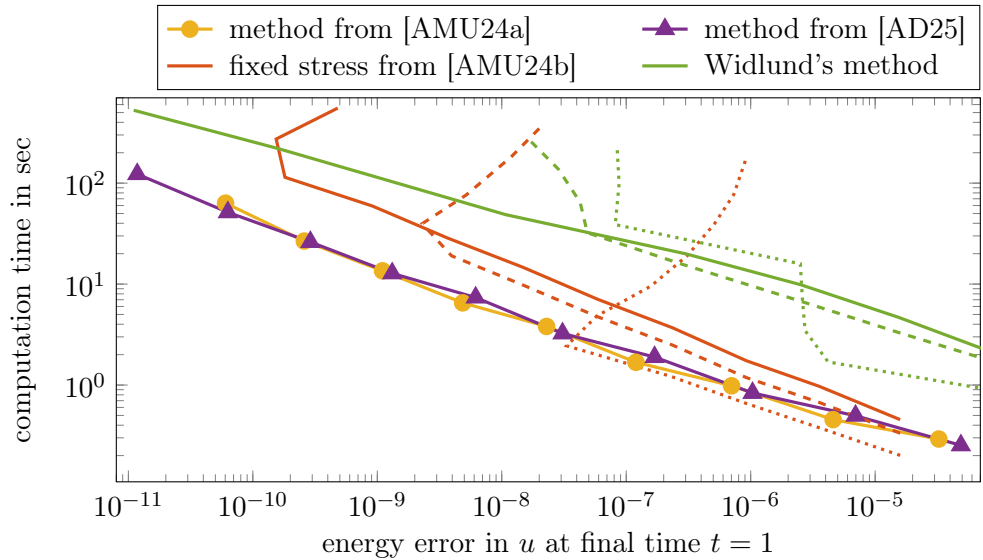

In this second experiment, we simulate the equations of poroelasticity on the time interval $[0, 1]$ with right-hand sides $f(t) = 1$, $g(t) = \sin(\pi t)$. Apart from that, we use the same physical parameters and initial data as before. The computation is done on a slightly coarser spatial mesh, leading to the dimensions $n = 7.938$ and $m = 3.969$. At this point, we would like to emphasize that all subsystems are solved with Matlab's backslash operator. The resulting computation times are shown in Figure~\ref{fig:poro:computationTimes}. Overall, it shows that all methods behave similarly as long as sufficient inner iteration steps are performed. If the parameter $K$ is too small compared to the time step size, then the fixed stress scheme even diverges. The same can be observed for Widlund's method. Here, we have set $\tol = c\,\tau^2$ for three different values of $c$ such that the number of inner iterations automatically increases for smaller time step sizes. 
Generally, one can say that the methods based on BDF-2 are slightly faster at the price that these methods do not guarantee the dissipation inequality on the discrete level. 
%

%
\section{Conclusion}
\label{sec:conclusion}
We have shown that the methods of Widlund and Rapoport can be employed to linear systems arising from the midpoint rule applied to general energy-based systems, if the systems are slightly modified and rescaled. With this, it is possible to derive an efficient implementation of a dissipation-preserving time stepping scheme of second order. Numerical examples examine different modeling strategies and compare the performance of the introduced methods with GMRES. 
In the case of poroelasticity, where the system matrix decouples due to its block structure, we show that the method of Widlund is compatible with well-established (but possibly non-dissipative) iterative schemes. 

Possible directions for future research are the extension to infinite-dimensional or to nonlinear energy-based DAE systems.
%
%
\section*{Acknowledgments} 
RA acknowledges support by the Deutsche Forschungsgemeinschaft (DFG, German Research Foundation) - 467107679.  
Major parts of this work were carried out while PS was affiliated with the University of Potsdam and supported by the DFG within the Sonderforschungsbereich 1294 “Data Assimilation -- The Seamless Integration of Data and Models” (Project 318763901).

Finally, the authors would like to thank Jörg Liesen for giving valuable comments on the paper.
%
%
\bibliographystyle{alpha}
\bibliography{references}
%
%
\appendix
\section{Implementation of Widlund's and Rapoport's Methods}
\label{app:Implementation}
In this section, we discuss some details of the implementation of the methods of Widlund and Rapoport, following \cite{SzyW93,Wid78} and  \cite{Rap78}, respectively. Note that, in contrast to these works, we consider $\Amat = \Hmat+\Smat$ rather than $\Amat = \Hmat-\Smat$.

As already mentioned, the methods are described by the projection processes \eqref{eq:Widlund} and \eqref{eq:Rapoport}, respectively. 
Since both methods revolve around the Krylov subspace $\calK_k(\Kmat,\widehat{r}_0)$, we  use the Arnoldi algorithm to compute an $\Hmat$-orthogonal basis of this Krylov space. To be more precise, we compute the basis vectors $v_1,\dots,v_k$ satisfying 	
\begin{equation*}
	\Kmat \Vmat_k 
	= \Vmat_k \Tmat_k + t_{k+1,k}v_{k+1}e_k^T 
	= \Vmat_{k+1}\Tmat_{k,k+1} 
	\quad \text{and}\quad 
	\Vmat_k^T \Hmat \Vmat_k 
	= \Dmat_k
	= \diag(\rho_1,\dots,\rho_n).
\end{equation*}
Here, $\Tmat_k\in\R^{k,k}$ is an unreduced upper Hessenberg matrix, $\Tmat_{k,k+1}=[\Tmat_k^T\; t_{k+1,k}e_k ]^T\in\R^{k+1,k}$, $\Vmat_\ell = [v_1,\dots,v_\ell]$ for $\ell=k,k+1$, and $\rho_1,\dots,\rho_n > 0$ are real and positive values. 
Further, $v_{k+1}$ is $\Hmat$-orthogonal to $\calK(\Kmat,\widehat{r}_0)$, i.e., $\Vmat_k^T\Hmat v_{k+1}=0$. 
This orthogonal basis is computed by a three-term recurrence, which can be seen from 
\begin{equation*}
	\Vmat_k^T\Smat\Vmat_k 
	= \Vmat_k^T \Hmat \Kmat \Vmat_k 
	= \Vmat_k^T\Hmat \Vmat_{k+1}\Tmat_{k,k+1}
	= \Dmat_k \big[\, \Imat_k\ 0\, \big]\ \Tmat_{k,k+1} 
	= \Dmat_k\Tmat_k.
\end{equation*}
Since $\Vmat_k^T\Smat\Vmat_k\in\R^{k,k}$ is skew-symmetric, so is $\Dmat_k\Tmat_k$ and, hence, $\Tmat_k$ has to be tridiagonal. More precisely, we have 
\begin{equation*}
	\Tmat_k 
	= \tri(\alpha_j,0,\beta_j) 
	= \begin{bmatrix}
		0 & \beta_2 \\ 
		\alpha_1 & 0 & \ddots \\ 
		&\ddots &\ddots&\ddots \\ 
		&&\ddots & 0 & \beta_{k} \\ 
		&&&\alpha_{k-1} & 0
	\end{bmatrix}
\end{equation*} 
with $\alpha_j,\beta_{j+1}\in\R$ for $j=1,\dots,k-1$. 
In particular, the recursion is given by 		
\begin{equation}\label{eq:Lanzsos-rec}
	\Kmat v_j 
	= \beta_{j} v_{j-1} + \alpha_j v_{j+1} 
	\iff 
	\alpha_jv_{j+1} 
	= \Kmat v_j - \beta_j v_{j-1}
\end{equation}
for $j=1,\dots,k$, where $v_{0}=0$. Note further that $\alpha_{j}\neq0$ for all $j=1,\dots,k-1$, since the matrix $\Tmat_k$ is unreduced. 
Due to the $\Hmat$-orthogonality of the vectors $v_j$, we see that 
\begin{align*}
	\rho_{j+1}\alpha_j  
	&= v_{j+1}^T \Hmat \Kmat v_j
	= v_{j+1}^T\Smat v_j \qquad \text{for }j = 1,\dots,k-1, \\ 
	\rho_{j-1}\beta_{j} 
	&= v_{j-1}^T \Hmat \Kmat v_j
	= v_{j-1}^T\Smat v_j \qquad \text{for }j = 2,\dots,k.
\end{align*}
This implies $\rho_j\beta_{j+1} = -\rho_{j+1}{\alpha_j}$ and, hence, $\beta_{j+1} = - {\alpha_j}\frac{\rho_{j+1}}{\rho_j}$ for $j=1,\dots,k$. 

The key observation here is that the scaling parameters $\alpha_j$ can be chosen arbitrarily (different choices just lead to different values of $\rho_j$). In order to derive the methods of Widlund and Rapoport, we use two different choices of these parameters. 
%
%
\subsection{Implementation of Widlund's method}

In order to implement Widlund's method, we choose $\alpha_1 = -1$ and 		
\begin{equation*}
	\alpha_{j+1} 
	= -1 + \alpha_j\frac{\rho_{j+1}}{\rho_j} 
	\qquad \text{for } j=1,\dots,k-2.
\end{equation*}
This way, we have 		
\begin{equation*}
	\alpha_j + \beta_j = -1 +  \alpha_{j-1}\frac{\rho_{j}}{\rho_{j-1}} -  \alpha_{j-1}\frac{\rho_{j}}{\rho_{j-1}} = -1 
	\qquad \text{for } j = 2,\dots k-1.
\end{equation*}

In order to compute the iterates of the linear system, we note that by \eqref{eq:Widlund}, we have $x_k = x_0+\Vmat_ky_k$ for some $y_k\in\R^k$ and that 		
\begin{align*}
	0 
	= \Vmat_k^T r_k 
	= \Vmat_k^T \Hmat\Hmat^{-1}(b - \Amat x_k) 
	= \Vmat_k^T \Hmat\Hmat^{-1}(r_0 - \Amat \Vmat_ky_k) 
	= \Vmat_k^T \Hmat\widehat{r}_0 - \Dmat_k(\Imat_n+\Tmat_k)y_k.
\end{align*}
Hence, in order to determine $y_k$, we have to solve the system 		
\begin{equation*}
	(\Imat_k+\Tmat_k)y_k 
	= e_1,
\end{equation*}
where we used that $\Dmat_k^{-1} \Vmat_k^T\Hmat\widehat{r}_0 = e_1$. 
This system can be solved by considering the $LU$-decomposition 		
\begin{equation*}
	\Imat_k+\Tmat_k 
	= 
	\begin{bmatrix}
		1 \\ -1 & 1 \\ &\ddots & \ddots \\ &&-1 & 1
	\end{bmatrix}
	\begin{bmatrix}
		-\alpha_1 & \beta_2 \\ &\ddots & \ddots \\ && -\alpha_{k-1} & \beta_{k} \\ &&& 1+\beta_{k}
	\end{bmatrix}
	\eqqcolon \Lmat_k \Umat_k.
\end{equation*}
By forward elimination, one directly sees that $w_k = [1,\, \dots,\, 1]^T \in \R^k$ solves $\Lmat_k w_k = e_1$. In order to solve the system $\Umat_ky_k = w_k$ via  backward substitution,  we note that $\alpha_{k-1}=-1-\beta_{k-1}$ and, hence,  
\begin{equation*}
	\Umat_k = 
	\begin{bmatrix}
		\Umat_{k-1} & \beta_{k}e_{k-1} \\ 0 & 1+\beta_{k}
	\end{bmatrix}.
\end{equation*}
In particular, we have 		
\begin{equation*}
	\Umat_k^{-1} = 
	\begin{bmatrix}
		\Umat_{k-1}^{-1} & -\frac{\beta_k}{1+\beta_k}\Umat_{k-1}^{-1}e_{k-1} \\ 
		0 & \frac{1}{1+\beta_{k}}
	\end{bmatrix}.
\end{equation*}
\begin{remark}
Note that $\Dmat_k(\Imat_k+\Tmat_k) = \Dmat_k + \Dmat_k \Tmat_k$ is the sum a symmetric and a skew-symmetric matrix, where the symmetric part is positive definite. Thus $(\Imat_k+\Tmat_k)$ has to be nonsingular and in particular we need to have $\beta_k\neq -1$.
\end{remark}
We now get 	
\begin{equation*}
	y_k = \Umat_k^{-1} w_k = 
	\begin{bmatrix}		
		\Umat_{k-1}w_{k-1} - \frac{\beta_k}{1+\beta_k}\Umat_{k-1}^{-1}e_{k-1} \\ 
		\frac{1}{1+\beta_{k}}
	\end{bmatrix}
	= 
	\begin{bmatrix}		
		y_{k-1} - \frac{\beta_k}{1+\beta_k}\Umat_{k-1}^{-1}e_{k-1} \\ 
		\frac{1}{1+\beta_{k}}
	\end{bmatrix}
	= 
	\begin{bmatrix}
		y_{k-1}\\0
	\end{bmatrix}
	+ 
	u_k
\end{equation*}
with
\begin{equation*}
	u_k = \Umat_k^{-1}e_k =
	\begin{bmatrix}		
		-\frac{\beta_k}{1+\beta_k}\Umat_{k-1}^{-1}e_{k-1} \\ 
		\frac{1}{1+\beta_{k}}
	\end{bmatrix}
	= 
	\begin{bmatrix}		
		-\frac{\beta_k}{1+\beta_k}u_{k-1} \\ 
		\frac{1}{1+\beta_{k}}
	\end{bmatrix}.
\end{equation*}
In particular, we have 		
\begin{align*}
	x_k-x_{k-1} 
	&= \Vmat_ky_k-\Vmat_{k-1}y_{k-1} = \Vmat_ku_k \\
	&= -\frac{\beta_k}{1+\beta_k}\Vmat_{k-1}u_{k-1} + \frac{1}{1+\beta_k} v_k \\ 
	&= -\left(1-\frac{1}{1+\beta_k} \right) (x_{k-1}-x_{k-2}) + \frac{1}{1+\beta_k} v_k
\end{align*}
for $k\ge 2$. Thus, we obtain 
\begin{equation*}
	x_k 
	= x_{k-2} +\underbrace{\frac{1}{1+\beta_k}}_{\eqqcolon\omega_k} \big(x_{k-1}-x_{k-2}+v_k\big)
\end{equation*}
for $k\ge 2$. By setting $x_{-1}=0$ and $\omega_1 = 1$ ($\beta_1=0$), the equation also holds for $k=1$. This is due to the fact that the projection process implies $x_1\in x_0+\calK_1(\Kmat,\widehat{r}_0)$, i.e., $x_1=x_0+\mu  \widehat{r}_0$ and 
\begin{equation*}
	0
	= \widehat{r}_0^T(b-\Amat x_1) 
	= \widehat{r_0}^T(b-\Amat x_0 - \mu \Amat\widehat{r}_0) 
	= \| \widehat{r}_0\|_\Hmat^2 - \mu (\widehat{r}_0^T\Hmat \widehat{r}_0+ \widehat{r}_0^T\Smat\widehat{r}_0) 
	= \|\widehat{r}_0\|_\Hmat^2 - \mu \|\widehat{r}_0\|_\Hmat^2.
\end{equation*}

In order to compute $x_k$ efficiently, we need to determine computable formulas for  $\omega_k$ and $v_k$. Starting with $\omega_k$, we note that $\omega_1=1$ and 
\begin{equation*}
	\omega_k 
	= \frac{1}{1+\beta_k} 
	= \frac{1}{-\alpha_k} 
	= \frac{1}{-(-1+\alpha_{k-1} \frac{\rho_{k}}{\rho_{k-1}})} 
	= \left(1 -\alpha_{k-1} \frac{\rho_{k}}{\rho_{k-1}} \right)^{-1} 
	= \left(1 +\frac{1}{\omega_{k-1}} \frac{\rho_{k}}{\rho_{k-1}} \right)^{-1}
\end{equation*}
for $k\ge2$.
Regarding the $v_k$, we observe that $v_1=\widehat{r}_0=\Hmat^{-1}r_0$. Inductively, we show that $v_{k+1}=\Hmat^{-1}r_k$. First,
\begin{equation*}
	r_k 
	= b - \Amat x_k 
	= b - \Amat(x_{k-2} +\omega_k(x_{k-1}-x_{k-2}+v_k)) 
	= r_{k-2} -\omega_k(r_{k-2}-r_{k-1} + \Amat v_k).
\end{equation*}
Multiplying from the left by $\Hmat^{-1}$ and assuming that $v_{\ell+1}=\Hmat^{-1}r_\ell$ for all $\ell<k$, we get 		
\begin{align*}
	\Hmat^{-1}r_k 
	&= \omega_k\left( \big(\tfrac{1}{\omega_k}-1 \big)\Hmat^{-1}r_{k-2} + \Hmat^{-1}r_{k-1}- v_k - \Kmat v_k \right) \\
	&= \frac{1}{-\alpha_k} \big( -\Kmat v_k + (-1-\alpha_k)v_{k-1} + v_k -v_k  \big) \\ 
	&= \frac{1}{-\alpha_k} \big( -\Kmat v_k + \beta_k v_{k-1} \big)  
	= \frac{-\alpha_{k}v_{k+1}}{-\alpha_k}
	= v_{k+1}.
\end{align*}
In summary, we obtain \Cref{alg:widlund} presented in Section~\ref{sec:wid_rap:wid}.
%
%
\subsection{Implementation of Rapoport's method} \label{sec:rapoport}
In order to implement the method of Rapoport, we choose the parameters $\alpha_j$ such that $\rho_j=1$. Thus it holds $\Vmat_k^T\Hmat\Vmat_k=\Imat_k$ and         
\begin{equation*}
    \Tmat_k = 
    \begin{bmatrix}
        0 & -\alpha_1 \\ 
        \alpha_1 & 0 & \ddots \\ 
        & \ddots & \ddots & \ddots \\
        && \ddots & 0 & -\alpha_{k-1} \\ 
        &&&\alpha_{k-1} & 0
    \end{bmatrix}.
\end{equation*}
Now the minimal residual property \eqref{eq:Rapoport_min} leads to the least squares problem 
\begin{align*}
	\min_{z\in x_0+\calK_k(\Kmat,\widehat{r}_0)} \norm{b-\Amat z}_{\Hmat^{-1}}
	&= \min_{y\in \R^{k}} \norm{r_0 - \Amat \Vmat_ky}_{\Hmat^{-1}} \\
	&= \min_{y\in \R^{k}} \norm*{r_0 - \Hmat\Vmat_{k+1} 
	\begin{bmatrix}
		\Imat_k+\Tmat_k \\
		\alpha_k e_{k}^T
	\end{bmatrix}
	y}_{\Hmat^{-1}} \\
	&= \min_{y\in \R^{k}} \norm*{\norm{\widehat{r}_0}_{\Hmat} e_1 -  
	\begin{bmatrix}
		\Imat_k+\Tmat_k \\
		\alpha_k e_{k}^T
	\end{bmatrix}
	y}
\end{align*}
Note that the last term is a least squares problem in the Euclidean norm instead of the $\Hmat^{-1}$-norm. 
We follow the general steps outlined in \cite[Sect.~3]{Fre90} to solve this least squares problem by a $QR$-decomposition which is computed via a sequence of Givens rotations	
\begin{equation*}
	\Gmat_k = 
	\begin{bmatrix}
		\Imat_{k-1} & 0 & 0 \\ 0 & c_k & s_k \\ 0 &-s_k & c_k
	\end{bmatrix}
	\quad\text{with}\quad 
	c_k^2+s_k^2 = 1.
	\end{equation*}
Here, each iteration only results in the application of one additional Givens rotation, as 
\begin{equation*}
	\begin{bmatrix}
		\Imat_k+\Tmat_k \\ \alpha_k e_k^{T}
	\end{bmatrix}
	= 
	\begin{bmatrix}
		\Qmat_k & 0 \\ 0& 1
	\end{bmatrix}
	\begin{bmatrix}
		\Rmat_k \\ \alpha_k e_k^{T}
	\end{bmatrix} 
	= 
	\begin{bmatrix}
		\Qmat_k & 0 \\ 0& 1
	\end{bmatrix}
	\Gmat_k^T\Gmat_k
	\begin{bmatrix}
		\Rmat_k \\ \alpha_k e_k^{T}
	\end{bmatrix}
	=
	\Qmat_{k+1}
	\begin{bmatrix}
		\widetilde{\Rmat}_{k} \\ 0
	\end{bmatrix},
\end{equation*}
where $\Qmat_k$ and $\Rmat_k$ are  the $QR$ decomposition of $\Imat_k+\Tmat_k$ and $\Gmat_{k}$ is the Givens rotation that eliminates the entry in the $(k+1)$st row. 
In particular we have that 
\[
	\gamma_1 
	= \sqrt{\alpha_1^2+1}, \qquad 
	c_1 
	=\frac{1}{\gamma_1},  \qquad \text{and} \quad 
	s_1 = \frac{\alpha_1}{\gamma_1}
\]
and it follows inductively, that 
\[
	\gamma_{k+1} 
	= \sqrt{\alpha_{k+1}^2+\frac{1}{c_k^2}}, \qquad 
	c_{k+1} 
	= \frac{1}{c_k\gamma_{k+1}} 
	= \frac{c_{k-1}\gamma_k}{\gamma_{k+1}}, \qquad \text{and} \quad 
	s_{k+1} = \frac{\alpha_{k+1}}{\gamma_{k+1}}.
\] 
The $QR$ decomposition of the full matrix $\Imat_{k+1}+\Tmat_{k+1}$ is given by 		
\begin{equation*}
	\Imat_{k+1}+\Tmat_{k+1} = \Qmat_{k+1} 
	\begin{bmatrix}
		\widetilde{\Rmat}_{k} & r_{k+1} \\ 0&  c_{k-1}\gamma_{k}
	\end{bmatrix}
	= \Qmat_{k+1}\Rmat_{k+1}
	\;\text{with}\; 
	\begin{bmatrix}
		r_{k+1}\\  c_{k-1}\gamma_{k} 
	\end{bmatrix}
	= \Gmat_k\Gmat_{k-1}(-\alpha_{k}e_{k}+e_{k+1})
\end{equation*}
and it can be shown that 		
\begin{equation*}
	\widetilde{\Rmat}_{k} = 
	\begin{bmatrix}
		\gamma_1 & 0 & -\alpha_2s_1 \\ 
		& \ddots & \ddots & \ddots  \\
		&& \ddots &\ddots & -\alpha_{k-1}s_{k-2} \\ 
		&&&\ddots & 0 \\
		&&&&\gamma_{k}
	\end{bmatrix}
	\quad\text{and}\quad
	r_{k+1} = -\alpha_ks_{k-1}e_{k-1}.
\end{equation*}
From 		
\begin{equation*}
	\min_{y\in \R^{k}} \norm*{\norm{\widehat{r}_0}_{\Hmat} e_1 -  
	\begin{bmatrix}
		\Imat_k+\Tmat_k \\
		\alpha_k e_{k}^T
	\end{bmatrix}
	y}
	=
	\min_{y\in \R^{k}} \norm*{\norm{\widehat{r}_0}_{\Hmat} \Qmat_{k+1}^Te_1 -  
	\begin{bmatrix}
		\widetilde{\Rmat}_k \\
		0
	\end{bmatrix}
	y}
\end{equation*}
we see that the minimizer $y_k$ is given by solving the linear system 	
\begin{equation*}
	\widetilde{\Rmat}_{k} y_k = d_{k}
	\quad\text{with}\quad
	\begin{bmatrix}
	 	d_k \\ \delta_k
	 \end{bmatrix} 
	 = \norm{\widehat{r}_0}_{\Hmat} \Qmat_{k+1}^Te_1.
\end{equation*}
In particular, we can see that 		
\begin{equation*}
	\begin{bmatrix}
	 	d_k \\ \delta_k
	 \end{bmatrix}
	 = 
	 \Gmat_{k}\left(\norm{\widehat{r}_0}_{\Hmat}
	 \begin{bmatrix}
	 	\Qmat_k^T & 0\\ 0 &1
	 \end{bmatrix}
	 e_1\right) 
	 = 
	 \Gmat_{k}
	 \begin{bmatrix}
	 	d_{k-1} \\ \delta_{k-1} \\0
	 \end{bmatrix}
	 = 
	 \begin{bmatrix}
	 	d_{k-1} \\ c_k\delta_{k-1} \\-s_k\delta_{k-1}
	 \end{bmatrix}.
\end{equation*}
Let $z_k\in \R^{k}$, be the solution of $\widetilde{\Rmat}_k z_k =e_k.$
Since 		
\begin{equation*}
	\widetilde{\Rmat}_k e_{k} = \gamma_k e_k  -\alpha_{k-1}s_{k-2}e_{k-2},
\end{equation*}
we see that 		
\begin{equation*}
	z_k = \frac{1}{\gamma_k} \left(e_k +  \alpha_{k-1}s_{k-2} 
	\begin{bmatrix}
	 	z_{k-2} \\ 0\\ 0
	\end{bmatrix} 
	\right).
\end{equation*}
In particular, we obtain 		
\begin{equation*}
	d_k = 
	\begin{bmatrix}
		d_{k-1} \\ c_k\delta_{k-1}
	\end{bmatrix}
	=
	\widetilde{\Rmat}_k
	\left(
	\begin{bmatrix}
		y_{k-1} \\ 0
	\end{bmatrix}
	+ 
	c_k\delta_{k-1} z_k
	\right)
\end{equation*} 		
and, thus, 	
\begin{equation*}
	y_k =
	\begin{bmatrix}
		y_{k-1} \\ 0
	\end{bmatrix}
	+ 
	c_k\delta_{k-1} z_k.
\end{equation*}
Setting $p_k = \Vmat_{k}z_k$, we obtain the update formulas 		
\begin{equation*}
	p_{k+1} 
	= \frac{1}{\gamma_{k+1}} \left(v_{k+1} +  \alpha_{k}s_{k-1} 
	p_{k-1}
	\right)
\end{equation*}
and 		
\begin{equation*}
	x_{k+1} 
	= x_0 + \Vmat_{k+1}y_{k+1} = x_k + c_{k+1}\delta_k p_{k+1}
\end{equation*}

Finally, the $\Hmat^{-1}$ norm of the residual can be easily updated throughout the iterations, since  $\norm{b-\Amat x_{k+1}}_{\Hmat^{-1}} = \delta_{k+1} = -s_{k+1}\delta_k$.


\end{document}

%% file: def.tex
\definecolor{myBlue}{RGB}{30,144,255} 
\definecolor{myGreen}{RGB}{69,169,0} 
\definecolor{myRed}{RGB}{165,12,42} 
\definecolor{myOrange}{RGB}{225,92,22} 
\definecolor{mycolor0}{rgb}{0.12156862745098,0.466666666666667,0.705882352941177} 
\definecolor{mycolor1}{rgb}{0.00000,0.44700,0.74100}
\definecolor{mycolor4}{rgb}{0.85000,0.32500,0.09800}
\definecolor{mycolor5}{rgb}{0.49400,0.18400,0.55600}
\definecolor{mycolor2}{rgb}{0.92900,0.69400,0.12500}
\definecolor{mycolor3}{rgb}{0.46600,0.67400,0.18800}
\definecolor{mycolor6}{rgb}{0.30100,0.74500,0.93300}
\definecolor{mycolor7}{rgb}{0.63500,0.07800,0.18400}

\newcommand{\delete}[1]{ }

\def\R{\mathbb{R}}

\DeclareMathOperator{\sspan}{span}

\DeclareMathOperator{\diag}{diag}
\DeclareMathOperator{\tri}{tri}

\newcommand{\calK}{\ensuremath{\mathcal{K}} }

\newcommand{\f}{\ensuremath{f}}
\newcommand{\g}{\ensuremath{g}}
\newcommand{\Amat}{\mathbf{A}}
\newcommand{\Bmat}{\mathbf{B}}
\newcommand{\Cmat}{\mathbf{C}}
\newcommand{\Dmat}{\mathbf{D}}
\newcommand{\Emat}{\mathbf{E}}

\newcommand{\Gmat}{\mathbf{G}}
\newcommand{\Hmat}{\mathbf{H}}
\newcommand{\Imat}{\mathbf{I}}
\newcommand{\Jmat}{\mathbf{J}}
\newcommand{\Kmat}{\mathbf{K}}
\newcommand{\Lmat}{\mathbf{L}}
\newcommand{\Mmat}{\mathbf{M}}

\newcommand{\Qmat}{\mathbf{Q}}
\newcommand{\Rmat}{\mathbf{R}}
\newcommand{\Smat}{\mathbf{S}}
\newcommand{\Tmat}{\mathbf{T}}
\newcommand{\Umat}{\mathbf{U}}
\newcommand{\Vmat}{\mathbf{V}}

\newcommand{\Nullmat}{\mathbf{0}} 

\newcommand{\tol}{\operatorname{TOL}}

\newcommand{\hamiltonian}{\mathcal{H}}

\newcommand{\inputPH}{u}

\newcommand{\state}{z}

\DeclarePairedDelimiter{\norm}{\lVert}{\rVert}
\DeclarePairedDelimiter{\abs}{\lvert}{\rvert}

\DeclareFontFamily{U}{matha}{\hyphenchar\font45}
\DeclareFontShape{U}{matha}{m}{n}{
	<-6> matha5 <6-7> matha6 <7-8> matha7
	<8-9> matha8 <9-10> matha9
	<10-12> matha10 <12-> matha12
}{}
\DeclareSymbolFont{matha}{U}{matha}{m}{n}

\DeclareFontFamily{U}{mathx}{\hyphenchar\font45}
\DeclareFontShape{U}{mathx}{m}{n}{
	<-6> mathx5 <6-7> mathx6 <7-8> mathx7
	<8-9> mathx8 <9-10> mathx9
	<10-12> mathx10 <12-> mathx12
}{}
\DeclareSymbolFont{mathx}{U}{mathx}{m}{n}

\DeclareMathDelimiter{\vvvert} {0}{matha}{"7E}{mathx}{"17}%

%% file: pics/biharmonic_1e-2.tex
%
%
\begin{tikzpicture}

\begin{axis}[%
width=0.26\textwidth,
height=0.202\textwidth,
at={(0\textwidth,0\textwidth)},
scale only axis,
xmin=0,
xmax=10,
xlabel style={font=\color{white!15!black}},
xlabel={iterations},
ymode=log,
ymin=1e-06,
ymax=2,
yminorticks=true,
ylabel style={font=\color{white!15!black}},
ylabel={rel.~residual norm},
axis background/.style={fill=white},
title style={font=\bfseries},
title={$\tau = \nnu^{-1} = 10^{-2} $},
trim axis left,
trim axis right
]
\addplot [color=mycolor1, very thick, forget plot]
  table[row sep=crcr]{%
0	1\\
1	0.10131130431435\\
2	0.00389489338302265\\
3	3.16658310741282e-05\\
4	5.58632430050129e-05\\
4.93173936690438	2.3436729115921e-07\\
};
\addplot [color=mycolor2,dashed, very thick, forget plot]
  table[row sep=crcr]{%
0	1\\
1	0.673256773209364\\
2	0.00389489316591398\\
3	8.80484555002722e-05\\
4	5.58360094106299e-05\\
4.98497643143713	2.3436729115921e-07\\
};
\addplot [color=mycolor3,dash dot, very thick, forget plot]
  table[row sep=crcr]{%
0	1\\
1	0.997570327679764\\
2	0.00389489589357955\\
3	0.00322611520820477\\
4	6.17380466226343e-05\\
5	2.28413799899936e-05\\
5.89451215600595	2.3436729115921e-07\\
};
\addplot [color=mycolor4,dotted, very thick, forget plot]
  table[row sep=crcr]{%
0	1\\
1	0.893888979494531\\
2	0.667478175473968\\
3	0.475852104924866\\
4	0.352004511524465\\
5	0.27428614812223\\
6	0.224145567200038\\
7	0.190676218546143\\
8	0.167658481757429\\
9	0.151431434992792\\
10	0.139750784584662\\
11	0.13118952178486\\
12	0.124812739509371\\
13	0.119992955591026\\
14	0.116300516072252\\
16	0.111188175069149\\
18	0.107975641735815\\
21	0.105100056336567\\
25	0.103084950563578\\
31	0.101692258663981\\
41	0.100799293335505\\
};
\end{axis}
\end{tikzpicture}%

%% file: pics/biharmonic_1e-4.tex
%
%
%
\begin{tikzpicture}

\begin{axis}[%
width=0.26\textwidth,
height=0.202\textwidth,
at={(0\textwidth,0\textwidth)},
scale only axis,
xmin=0,
xmax=20,
xlabel style={font=\color{white!15!black}},
xlabel={iterations},
ymode=log,
ymin=1e-06,
ymax=2,
ytick={1e-06,1e-03},
yticklabels={\empty},
yminorticks=true,
axis background/.style={fill=white},
title style={font=\bfseries},
title={$\tau = \nnu^{-1} = 10^{-4} $},
trim axis left,
trim axis right
]
\addplot [color=mycolor1, very thick, forget plot]
  table[row sep=crcr]{%
0	1\\
1	0.101321129988216\\
2	0.00101616280190012\\
3	2.42505662686398e-07\\
};
\addplot [color=mycolor2, very thick, dashed, forget plot]
  table[row sep=crcr]{%
0	1\\
1	0.995153272893541\\
2	0.00101616280199085\\
3	0.000694864864539417\\
4	0.000694701031427356\\
5	0.000135133816939476\\
6	0.000135122755672316\\
7	8.01238292060788e-06\\
8	8.01207220673515e-06\\
9	1.96603982281345e-06\\
10	1.96603988937627e-06\\
11	8.57973869024128e-07\\
};
\addplot [color=mycolor3, very thick, dash dot, forget plot]
  table[row sep=crcr]{%
0	1\\
1	0.999999756477326\\
2	0.00101616280190578\\
3	0.00101614178014298\\
4	0.000724817686797841\\
5	0.000724690176532671\\
6	0.000163514898327128\\
7	0.000163402017851947\\
8	4.90831340586821e-05\\
9	1.24464552224287e-05\\
10	1.00064311249977e-05\\
11	9.86710500608308e-06\\
12	2.53625937086283e-07\\
};
\addplot [color=mycolor4, very thick, dotted, forget plot]
  table[row sep=crcr]{%
0	1\\
7	0.996929145137275\\
8	0.983283551010728\\
9	0.92708057171729\\
10	0.851546051860516\\
12	0.653422534653723\\
16	0.427737034944087\\
17	0.387635837794071\\
18	0.348529782457298\\
19	0.317189010199478\\
21	0.273596651507467\\
22	0.253125658348942\\
23	0.237151775809062\\
25	0.211211164036104\\
27	0.190544308433615\\
29	0.174723470629456\\
32	0.156841728700626\\
34	0.148195577845753\\
36	0.140846431621797\\
39	0.1324348714702\\
41	0.128091446920614\\
};
\end{axis}
\end{tikzpicture}%

%% file: pics/biharmonic_1e-6.tex
%
%
%
\begin{tikzpicture}

\begin{axis}[%
width=0.26\textwidth,
height=0.202\textwidth,
at={(0\textwidth,0\textwidth)},
scale only axis,
xmin=0,
xmax=40,
xlabel style={font=\color{white!15!black}},
xlabel={iterations},
ymode=log,
ymin=1e-06,
ymax=2,
ytick={1e-06,1e-03},
yticklabels={\empty},
yminorticks=true,
axis background/.style={fill=white},
title style={font=\bfseries},
title={$\tau = \nnu^{-1} = 10^{-6} $},
trim axis left,
trim axis right
]
\addplot [color=mycolor1, very thick, forget plot]
  table[row sep=crcr]{%
0	1\\
1	6.11744434995084e-05\\
2	0.999951115227899\\
2.83550155465824	2.3436729115921e-07\\
};
\addplot [color=mycolor2, very thick, dashed, forget plot]
  table[row sep=crcr]{%
0	1\\
2	0.999951115227901\\
2.99161617454174	2.3436729115921e-07\\
};
\addplot [color=mycolor3, very thick, dash dot, forget plot]
  table[row sep=crcr]{%
0	1\\
5	0.997919820384291\\
9	0.983150505448541\\
10	0.944057245789395\\
11	0.943111857155904\\
12	0.754095199704047\\
13	0.753550146697202\\
14	0.536767347870313\\
15	0.528020842487446\\
16	0.243044435925792\\
17	0.23574140034437\\
18	0.0535752525598507\\
19	0.0411256761735381\\
20	0.011235177265823\\
21	0.0112063018735523\\
22	0.00270272603041427\\
23	0.00173349201799847\\
24	0.000504513809132793\\
25	0.000240459357596501\\
26	6.18206435475171e-05\\
27	5.48700519453082e-05\\
28	9.14807007406043e-06\\
29	9.21084297308186e-06\\
30	1.08094718343914e-06\\
31	8.52465402276712e-07\\
};
\addplot [color=mycolor4, very thick, dotted, forget plot]
  table[row sep=crcr]{%
0	1\\
1	0.0840517199308957\\
2	0.0209677956339146\\
3	0.00536584481735382\\
4	0.00238449452615874\\
5	0.00130488631053739\\
6	0.000909737647069392\\
7	0.000728728766832995\\
8	0.000648262377812658\\
9	0.000618788146041789\\
10	0.000608886637322161\\
13	0.000604740917307618\\
23	0.000603776519011978\\
41	0.0006037676061559\\
};
\end{axis}
\end{tikzpicture}%

%% file: pics/biharmonic_widlund_norms.tex
%
%
\begin{tikzpicture}

\begin{axis}[%
width=0.28\textwidth,
height=0.25\textwidth,
at={(0\textwidth,0\textwidth)},
scale only axis,
xmin=0,
xmax=50,
xlabel={iterations},
ymode=log,
ymin=1e-12,
ymax=100000,
yminorticks=true, 
ylabel={rel.~resdiual norm},
ytick={1e-06,1e-12,1e-0,1e4},
axis background/.style={fill=white},
title={Widlund},
trim axis left,
trim axis right
]
\addplot [color=mycolor1, very thick, forget plot]
  table[row sep=crcr]{%
0	1\\
1	0.0514553251883739\\
2	0.86102986751899\\
3	9.82050455265596e-08\\
4	0.0145586917716256\\
5	9.36679141549007e-08\\
6	0.000199904329675308\\
7	1.05579431977034e-06\\
8	0.00019425058998578\\
9	7.0725159939721e-08\\
10	1.17550299217904e-05\\
11	6.93983289059086e-08\\
12	7.92680836306081e-06\\
13	2.99922142205928e-08\\
14	4.12570761149337e-07\\
15	3.35586145268125e-08\\
16	3.73902950729467e-07\\
17	6.26366349000208e-09\\
18	1.38457289913881e-08\\
19	6.99741778718775e-09\\
20	1.38138522018952e-08\\
21	3.68276815588537e-09\\
22	1.33726325984279e-08\\
23	2.03891335617196e-09\\
24	1.15332568237509e-08\\
25	1.26762798067889e-09\\
26	1.0780346550846e-08\\
27	1.18850083235518e-09\\
28	9.30096900440652e-09\\
29	5.55378806254315e-10\\
30	2.97442287015852e-10\\
31	5.55667971993525e-10\\
32	2.83816533794554e-10\\
33	4.21485793293302e-10\\
34	2.73195249203235e-10\\
35	1.28695157219316e-10\\
36	2.74023166589377e-10\\
37	6.13437765019325e-11\\
38	6.30273906821968e-11\\
39	6.71993996960283e-11\\
40	6.60705155688993e-11\\
41	7.11562030973988e-11\\
42	5.74795467481007e-11\\
43	7.49727786460513e-11\\
44	6.27692572971223e-11\\
45	7.91501176766131e-11\\
46	5.10504507023198e-11\\
47	5.64719193339409e-11\\
49	6.13332690522539e-11\\
50	5.61960643807096e-11\\
};
\addplot [color=mycolor2, dashed, very thick, forget plot]
  table[row sep=crcr]{%
0	1\\
1	45.3122108179784\\
2	2.70674863016103e-05\\
3	3.54978987104739e-06\\
4	1.74230535058529e-05\\
5	4.07059323733808e-07\\
6	7.35009943967671e-07\\
7	2.07624918186831e-06\\
8	2.56072143781686e-06\\
9	8.77499083436447e-08\\
10	5.3586298009818e-07\\
11	1.21038548374847e-07\\
12	6.0566210320768e-08\\
13	1.06584891548962e-08\\
14	1.43344163853475e-07\\
15	3.21422904807872e-08\\
16	5.08307310545386e-09\\
17	8.82409241034477e-10\\
18	9.40276413452368e-10\\
19	1.32693094845445e-08\\
20	9.88238701562623e-10\\
21	2.95650555466708e-09\\
22	7.75688391172225e-10\\
23	1.09436898530758e-09\\
24	4.55522749174894e-10\\
26	2.83917445750944e-09\\
27	7.27751105694234e-10\\
28	3.37195348395006e-10\\
29	3.73528736433554e-11\\
30	3.67538856218242e-10\\
31	1.58738662392271e-10\\
32	6.00059754749974e-11\\
33	1.75951094225753e-10\\
34	1.04674645576494e-11\\
35	2.32885901548475e-10\\
36	1.59005143073334e-11\\
37	2.86502551947268e-12\\
38	2.27491212479006e-11\\
39	1.12453032419708e-11\\
40	1.65220433556179e-11\\
41	3.11331264773086e-12\\
42	1.15880321946978e-11\\
43	6.49441872955145e-12\\
44	2.16408286843284e-11\\
45	1.62222500984011e-12\\
46	1.12075708329444e-12\\
47	7.64234188054717e-12\\
48	3.69586837298207e-12\\
49	1.32656964377653e-12\\
50	3.3962954036851e-12\\
};
\addplot [color=mycolor5, dashdotted, very thick, forget plot]
  table[row sep=crcr]{%
0	1.00000000087421\\
1	20264.2367282915\\
2	2.39325783935769e-07\\
3	0.000526014574146544\\
4	5.91051892233384e-06\\
5	3.61493775291416e-05\\
6	8.8288421795305e-05\\
7	1.91986953647666e-05\\
8	0.00110293506309741\\
9	1.20481825519654e-05\\
10	1.05896757130071e-05\\
11	2.2571356946564e-05\\
12	1.45170389938765e-07\\
13	5.17755306467681e-07\\
14	6.38583153333601e-05\\
15	1.20597770156682e-06\\
16	7.98547963043647e-08\\
17	3.91465442839493e-08\\
18	2.60007968140392e-07\\
19	2.9655662103082e-06\\
20	2.95080490073419e-07\\
21	2.89029159078411e-07\\
22	1.2975086577212e-07\\
23	8.67496155110799e-08\\
24	6.01727416007894e-08\\
25	1.13600856050469e-07\\
26	1.22383273552403e-06\\
27	6.08765851329006e-08\\
28	2.32542438540263e-08\\
29	1.68454316802817e-09\\
30	1.61079702088115e-07\\
31	8.02478778728405e-09\\
32	8.97142813120225e-09\\
33	3.12554428028605e-08\\
34	2.36797583813516e-09\\
35	4.10720850836708e-08\\
36	5.45629022642776e-09\\
37	3.73253604406239e-10\\
38	5.61207414096011e-09\\
39	3.59844683809694e-09\\
40	3.2450914509405e-09\\
41	6.31879107671895e-10\\
42	7.40302244499937e-10\\
43	2.81495424374347e-09\\
44	1.17337018964609e-09\\
45	2.13901408277654e-10\\
46	7.92202719969086e-11\\
47	3.277477401905e-09\\
48	3.96344316118729e-10\\
49	2.60537948724886e-10\\
50	2.05961338028853e-10\\
};
\addplot [color=gray, dotted, thick, forget plot]
  table[row sep=crcr]{%
0	1\\
50	1\\
};
\end{axis}
\end{tikzpicture}%

%% file: pics/biharmonic_rapoport_norms.tex
%
%
\begin{tikzpicture}

\begin{axis}[%
width=0.28\textwidth,
height=0.25\textwidth,
at={(0\textwidth,0\textwidth)},
scale only axis,
xmin=0,
xmax=50,
xlabel={iterations},
ymode=log,
ymin=1e-12,
ymax=100000,
yticklabels={\empty},
yminorticks=true,
ytick={1e-06,1e-12,1e-0,1e4},
axis background/.style={fill=white},
title={Rapoport},
trim axis left,
trim axis right
]
\addplot [color=mycolor1, very thick, forget plot]
  table[row sep=crcr]{%
0	1\\
1	0.999513191956026\\
2	0.861029867519117\\
3	0.0145586680025867\\
4	0.0145586680383882\\
5	0.000199872807224045\\
6	0.000199873306860863\\
7	1.22092692356568e-05\\
8	1.22092787201461e-05\\
9	7.93443836566351e-06\\
12	7.9007001576163e-06\\
13	3.73450509619598e-07\\
14	3.73901577049936e-07\\
15	1.36197819470614e-08\\
16	1.36815746613815e-08\\
17	1.07504652490793e-08\\
18	1.07620512255056e-08\\
19	3.08874716189264e-10\\
20	3.13708202049379e-10\\
22	3.11249652061913e-10\\
23	1.44402390839415e-10\\
24	1.44836361049359e-10\\
25	1.47503581265259e-10\\
26	1.48014191762136e-10\\
27	1.50690035712469e-10\\
28	1.50470781287276e-10\\
29	1.53137727945952e-10\\
30	1.53200458008889e-10\\
31	1.55602558800084e-10\\
34	1.55508345866702e-10\\
35	1.58785335594435e-10\\
50	1.58949874985706e-10\\
};
\addplot [color=mycolor2, dashed, very thick, forget plot]
  table[row sep=crcr]{%
0	1\\
1	0.999756566190664\\
2	2.70674862923977e-05\\
3	3.51965072038472e-06\\
4	3.44996069298626e-06\\
5	4.04230136141736e-07\\
6	3.51041758034771e-07\\
7	8.67629783443755e-08\\
8	8.57665144948208e-08\\
9	6.91378725790557e-08\\
10	4.86949241591843e-08\\
11	4.858710471205e-08\\
12	4.54505189351009e-08\\
13	9.88388747240791e-09\\
14	4.52105592594082e-09\\
15	8.61895167784746e-10\\
16	4.20345194926549e-10\\
17	3.73090234627608e-10\\
18	2.33633817963682e-10\\
19	4.30739451203033e-11\\
20	2.15305748806933e-11\\
21	2.24166623720022e-11\\
22	2.23735451670356e-11\\
23	1.9431150889593e-11\\
27	1.9425267078441e-11\\
38	1.94682409309895e-11\\
50	1.94682462694083e-11\\
};
\addplot [color=mycolor5, dashdotted, very thick, forget plot]
  table[row sep=crcr]{%
0	1.00000000087421\\
1	9.9153064391657\\
2	2.39802119198249e-07\\
3	0.000517120433209657\\
4	0.000496844978054084\\
5	3.52952439625685e-05\\
6	2.66180607852236e-05\\
7	1.19293994809318e-05\\
8	1.16578129625848e-05\\
9	4.19839577513294e-06\\
10	2.08295607196665e-06\\
11	2.42917646242622e-06\\
12	2.12510165633747e-06\\
13	4.64875083852753e-07\\
14	9.76789785275018e-08\\
15	3.73559596563916e-08\\
16	1.17127349021142e-08\\
17	1.95296947042489e-08\\
18	1.05955314131026e-08\\
19	8.06157478464975e-09\\
20	7.96159843733775e-09\\
21	8.4686612854167e-09\\
22	8.45784599304771e-09\\
23	7.9693125183335e-09\\
36	7.97115495181562e-09\\
50	7.96388033514549e-09\\
};
\addplot [color=gray, dotted,thick, forget plot]
  table[row sep=crcr]{%
0	1\\
50	1\\
};
\end{axis}
\end{tikzpicture}%

%% file: pics/biharmonic_pgmres_norms.tex
%
%
\begin{tikzpicture}

\begin{axis}[%
width=0.28\textwidth,
height=0.25\textwidth,
at={(0\textwidth,0\textwidth)},
scale only axis,
xmin=0,
xmax=50, 
xlabel={iterations},
ymode=log,
ymin=1e-12,
ymax=100000,
ytick={1e-06,1e-12,1e-0,1e4},
yticklabels={\empty},
yminorticks=true,
axis background/.style={fill=white},
title={prec.~GMRES},
trim axis left,
trim axis right
]
\addplot [color=mycolor1, very thick, forget plot]
  table[row sep=crcr]{%
0	1\\
1	0.999999997564769\\
2	0.861029867519512\\
3	0.861029585580849\\
4	0.828128924265229\\
5	0.828126218917224\\
6	0.784170255996221\\
8	0.620051698267057\\
9	0.618419280716413\\
10	0.257559540275708\\
11	0.257557519572609\\
12	0.143383043523891\\
13	0.137171241284151\\
14	0.0435824244871361\\
15	0.0273460993107515\\
16	0.00496043637402813\\
17	0.00495409140397666\\
18	0.00233538779473305\\
19	0.00231211410578569\\
20	0.000130608207533061\\
21	0.000123408155875521\\
22	1.23382291012952e-05\\
23	1.22817313775387e-05\\
24	1.17247759883582e-06\\
25	1.08882903811224e-06\\
26	1.17208506302913e-07\\
27	1.18149744674648e-07\\
29	1.16127340238124e-07\\
30	1.51599978914396e-07\\
31	2.5137150000159e-07\\
32	1.57718842110862e-07\\
33	1.21558750247643e-07\\
34	1.28269549173707e-07\\
35	1.46927797093714e-07\\
36	2.04756836154872e-07\\
37	1.18006588634503e-07\\
38	1.99970799694308e-07\\
39	1.44608694439007e-07\\
40	1.15934232319666e-07\\
41	1.9605211373983e-07\\
42	2.08371072505973e-07\\
43	2.58897105253983e-07\\
44	2.42340507304946e-07\\
45	1.20843448468971e-07\\
46	1.52241036225272e-07\\
47	3.14472118402979e-07\\
48	1.36154714121198e-07\\
49	1.17892829298447e-07\\
50	1.67860596722096e-07\\
};
\addplot [color=mycolor2, dashed, very thick, forget plot]
  table[row sep=crcr]{%
0	1\\
1	0.999999997564769\\
2	2.706748629171e-05\\
3	2.70674774286322e-05\\
4	2.60242302703107e-05\\
5	2.60241452538469e-05\\
6	2.4642509730514e-05\\
8	1.94845634970192e-05\\
9	1.9433265640756e-05\\
10	8.09350633053075e-06\\
11	8.09344283233186e-06\\
12	4.50563931725554e-06\\
13	4.31044042440833e-06\\
14	1.36952290658739e-06\\
15	8.59316639373243e-07\\
16	1.55875315596886e-07\\
17	1.55675933219212e-07\\
18	7.33865508754577e-08\\
19	7.26552049835376e-08\\
20	4.10419431645433e-09\\
21	3.87794210690726e-09\\
22	3.87723661298399e-10\\
23	3.85948340330985e-10\\
24	3.69533085074564e-11\\
25	3.43329633904837e-11\\
26	4.64874374544929e-12\\
28	4.64580435276414e-12\\
29	4.62195742461878e-12\\
30	5.54413336856727e-12\\
31	8.39239446314337e-12\\
32	5.71047612572664e-12\\
33	4.75824996295974e-12\\
34	4.9287184447521e-12\\
35	5.41959696725254e-12\\
36	7.03199430921934e-12\\
37	4.66921723990574e-12\\
38	6.89461749422372e-12\\
39	5.35690847503647e-12\\
40	4.61706393685484e-12\\
41	6.78190525106029e-12\\
42	7.13527509903952e-12\\
43	8.61540720930571e-12\\
44	8.12585748584867e-12\\
45	4.73997327039747e-12\\
46	5.561167703506e-12\\
47	1.02814511225424e-11\\
48	5.13367290125102e-12\\
49	4.66604814998805e-12\\
50	5.98994605068503e-12\\
};
\addplot [color=mycolor5, dashdotted, very thick, forget plot]
  table[row sep=crcr]{%
0	1.00000000087421\\
1	0.999999999656586\\
2	2.39323993798941e-07\\
3	2.39323928716317e-07\\
4	4.71281464252579e-08\\
5	4.71278489489107e-08\\
6	4.13898929065012e-08\\
7	2.87353590341155e-08\\
8	1.18271510752012e-08\\
9	1.18012318330117e-08\\
10	2.23160904602072e-09\\
11	2.23164356094383e-09\\
12	1.38605131185307e-09\\
13	1.33702969899907e-09\\
14	2.06200661791861e-10\\
15	1.5674530898114e-10\\
16	4.15971362286592e-11\\
17	4.73936771953655e-11\\
18	4.29180930769777e-11\\
19	2.84685593583846e-11\\
20	2.89915897681529e-11\\
21	2.40228276067811e-11\\
22	3.33978553710922e-11\\
23	3.27444227811621e-11\\
24	5.03702162094273e-11\\
25	3.09434742915106e-11\\
26	3.00999930102549e-11\\
27	4.01238473556922e-11\\
28	4.23659814261841e-11\\
29	3.54712280508413e-11\\
30	3.63088254672194e-11\\
31	3.39652949485702e-11\\
32	4.30595503176274e-11\\
33	3.97185749448114e-11\\
34	2.65285975751718e-11\\
35	4.51184515017641e-11\\
36	2.57109922538304e-11\\
37	4.17047022565902e-11\\
38	2.37824545414653e-11\\
39	2.79039278466798e-11\\
40	3.04597498896857e-11\\
41	3.43309012315e-11\\
42	2.91315844569238e-11\\
43	3.65037325890672e-11\\
44	3.31997367962651e-11\\
45	3.43672461828329e-11\\
46	2.70276391942172e-11\\
47	3.42892755537608e-11\\
48	3.52226296468775e-11\\
50	4.40411698327581e-11\\
};
\addplot [color=gray, dotted,thick , forget plot]
  table[row sep=crcr]{%
0	1\\
50	1\\
};
\end{axis}
\end{tikzpicture}%

%% file: pics/square_oneStep_nRef8.tex
%
%
\begin{tikzpicture}

\begin{axis}[%
width=4.5in,
height=1.7in,
scale only axis,
xmode=log,
xmin=2e-10,
xmax=0.01,
xminorticks=true,
xlabel={combined energy error in $u$ and $p$},
ymin=0.5,
ymax=8.3,
yminorticks=true,
ylabel={computation time in sec},
legend cell align={left},
axis background/.style={fill=white}
]
\addplot [color=mycolor1, very thick, mark=*, mark size=2.7]
  table[row sep=crcr]{%
0.00743019464748559	1.492923\\
0.000134598799561962	2.671299\\
1.648632119737e-06	3.490796\\
1.77787855591196e-08	4.723859\\
2.92385642813045e-10	5.718022\\
};
\addlegendentry{Widlund (backslash)}

\addplot [color=mycolor2, very thick, mark=triangle*, mark size=3.2]
table[row sep=crcr]{%
0.00747005997588053	2.110036\\
0.000139003581694943	3.517694\\
2.12165782970824e-06	4.66512\\
9.014309801563e-08	5.909216\\
1.59688781252073e-09	7.577368\\
4.29749195283855e-11	8.578843\\
};
\addlegendentry{Widlund (CG)}

\addplot [color=mycolor4, very thick, mark=square*, mark size=2.7]
table[row sep=crcr]{%
0.00743019464748637	1.025681\\
0.00743019464748637	1.011258\\
0.000134598799562056	1.167529\\
1.64863211972535e-06	1.255138\\
1.77787856023182e-08	1.497911\\
2.92385229147812e-10	1.570284\\
};
\addlegendentry{Widlund (Cholesky)}

\addplot [color=mycolor3, very thick, dashdotted]
table[row sep=crcr]{%
0.00590380641194635	2.456777\\
0.000816267341887614	2.892268\\
0.000120765283699966	3.361404\\
1.90791348384272e-05	4.025411\\
3.75833162747606e-06	4.438005\\
1.21465828997984e-07	5.4837\\
1.39846961304098e-08	5.860507\\
1.72083288384578e-09	6.620096\\
2.74162430111309e-10	7.005372\\
};
\addlegendentry{prec.~GMRES}

\addplot [color=gray, dashed]
  table[row sep=crcr]{%
0.1	3.154473\\
1.e-11	3.154473\\
};

\end{axis}

\end{tikzpicture}%

%% file: pics/square_times_nRef6_T1_om05.tex
%
%
%
\begin{tikzpicture}

\begin{axis}[%
width=4.5in,
height=2.0in,
scale only axis,
xmode=log,
xmin=8e-12,
xmax=7e-5,
xminorticks=true,
xlabel={energy error in $u$ at final time $t=1$},
ymode=log,
ymin=0.11,
ymax=700,
yminorticks=true,
ylabel={computation time in sec},
axis background/.style={fill=white},
legend columns = 2, 
legend style={at={(1.0,1.02)}, anchor=south east, legend cell align=left, align=left}
]

\addplot[color=mycolor2, very thick, mark=*, mark size=2.7]
  table[row sep=crcr]{%
3.24800130214256e-05	0.291903\\
4.59628004466635e-06	0.454486\\
7.02538274029882e-07	0.982828\\
1.19660868061548e-07	1.680882\\
2.29183873289841e-08	3.805122\\
4.85367054164595e-09	6.515844\\
1.10168508797689e-09	13.525277\\
2.59096759772038e-10	26.678874\\
6.03652605132591e-11	63.333028\\
};
\addlegendentry{\ method from \cite{AltMU24}}

\addplot[color=mycolor5, very thick, mark=triangle*, mark size=3.3]
  table[row sep=crcr]{%
4.88496675260677e-05	0.252766\\
6.95587461847436e-06	0.496249\\
1.03494602172384e-06	0.834754\\
1.68973909784671e-07	1.892521\\
3.07682773250242e-08	3.243809\\
6.16668477047102e-09	7.356068\\
1.31748235100517e-09	12.77156\\
2.89852564840269e-10	26.111715\\
6.29624228775902e-11	51.435795\\
1.17860379692282e-11	122.792929\\
};
\addlegendentry{\ method from \cite{AltD25}} 

\addplot [color=mycolor4, very thick, dotted, forget plot]
  table[row sep=crcr]{%
1.58893081903229e-05	0.201016\\
3.53555876959154e-06	0.378483\\
9.26158949776589e-07	0.655802\\
2.26555001173372e-07	1.222752\\
3.32057663255956e-08	2.446079\\
5.80746464214801e-08	4.82611\\
1.58385832579279e-07	9.582196\\
3.09095699607313e-07	19.108549\\
5.10357562445884e-07	39.3452\\
7.30901616705961e-07	77.030543\\
9.32845218816408e-07	185.080679\\
};

\addplot [color=mycolor4, very thick, dashed, forget plot]
  table[row sep=crcr]{%
1.58894625655438e-05	0.333118\\
3.53622418616281e-06	0.673851\\
9.28663844247205e-07	1.181857\\
2.35180771079709e-07	2.488028\\
5.90235250906554e-08	4.741579\\
1.4805385192893e-08	9.780673\\
3.94992707585698e-09	19.054128\\
2.29488420168968e-09	39.66554\\
5.16984154510003e-09	76.169536\\
1.19771416332423e-08	183.415052\\
2.09924884053642e-08	369.35196\\
};

\addplot [color=mycolor4, very thick]
  table[row sep=crcr]{%
1.58894625655344e-05	0.454274\\
3.53622418584532e-06	0.973289\\
9.28663835827598e-07	1.73591\\
2.35180588309676e-07	3.693098\\
5.90205368232581e-08	7.066035\\
1.47715607382345e-08	14.683771\\
3.69156142331502e-09	28.496197\\
9.14300892345392e-10	59.28216\\
1.81290007133308e-10	114.186096\\
1.53953880484397e-10	273.307139\\
4.80949521695891e-10	553.56224\\
};
\addlegendentry{\ fixed stress from \cite{AltMU24b}\qquad} 

\addplot[color=mycolor3, very thick, dotted, forget plot]
  table[row sep=crcr]{%
6.71546402783979e-05	0.939731\\
4.1890632598248e-06	1.695173\\
2.73830931421674e-06	3.459486\\
2.60711487063333e-06	6.8978\\
2.50567787992158e-06	15.872069\\
8.1634093707096e-08	38.851831\\
8.76649308710118e-08	98.515721\\
8.46714050625458e-08	213.23736\\
};

\addplot[color=mycolor3, very thick, dashed, forget plot]
  table[row sep=crcr]{%
8.66875118592525e-05	1.723235\\
1.53935391970749e-05	3.314628\\
2.42151861052488e-06	6.874957\\
2.79748738054418e-07	15.751563\\
4.871939287265e-08	32.229787\\
4.2607300201482e-08	64.698035\\
2.85823769084446e-08	136.997539\\
1.62524892530108e-08	276.314172\\
};

\addplot[color=mycolor3, very thick]
  table[row sep=crcr]{%
8.67564388066301e-05	2.113862\\
1.54979944105213e-05	4.672637\\
2.53651268215322e-06	9.833375\\
3.02148058937456e-07	19.91857\\
1.07243538387841e-08	48.76516\\
1.59149281550838e-10	220.196232\\
1.10394120088517e-11	526.912158\\
};
\addlegendentry{\ Widlund's method} 

\end{axis}

\end{tikzpicture}%